\documentclass[11pt]{article}

\usepackage[margin=1in]{geometry}
\usepackage[T1]{fontenc}
\usepackage{lmodern}
\usepackage{microtype}
\usepackage{setspace}
\usepackage{amsmath, amsthm, amssymb, amsfonts}
\usepackage{mathrsfs}
\usepackage{mathtools}
\usepackage{aliascnt}
\renewcommand{\arraystretch}{1.15}

\usepackage{enumitem}
\usepackage{needspace}
\usepackage{longtable}
\usepackage{tikz-cd}
\usepackage{extarrows}

\usepackage{booktabs}
\usepackage{array}
\usepackage{caption}
\usepackage{ragged2e}
\usepackage{adjustbox}

\newcolumntype{T}[1]{>{\RaggedRight\arraybackslash}p{#1}}
\newcolumntype{Q}[1]{>{\centering\arraybackslash}p{#1}}

\usepackage[backend=biber, bibstyle=trad-abbrv, citestyle=numeric,
  sorting=nty, giveninits=true]{biblatex}
\DeclareFieldFormat{titlecase}{#1}
\usepackage[colorlinks=true, citecolor=blue, linkcolor=blue, urlcolor=blue]{hyperref} 
\usepackage[nameinlink,capitalise,noabbrev]{cleveref}
\usepackage{comment}
\usepackage{makecell}

\newcommand{\cC}{\mathcal{C}}

\newcommand{\F}{\mathbb{F}}
\newcommand{\PP}{\mathbb{P}} 

\DeclareMathOperator{\wt}{wt}
\DeclareMathOperator{\Sym}{Sym}
\DeclareMathOperator{\Aut}{Aut}
\DeclareMathOperator{\GL}{GL}

\DeclareMathOperator{\Irr}{Irr}
\DeclareMathOperator{\Gal}{Gal}
\DeclareMathOperator{\nhyp}{nhyp}
\DeclareMathOperator{\PGL}{PGL}
\DeclareMathOperator{\Pic}{Pic}
\DeclareMathOperator{\Gr}{Gr}
\DeclareMathOperator{\Frob}{Frob}

\theoremstyle{plain}
\newtheorem{thm}{Theorem}[section]
\newaliascnt{prop}{thm}
\newtheorem{prop}[prop]{Proposition}
\aliascntresetthe{prop}
\newaliascnt{lemm}{thm}
\newtheorem{lemm}[lemm]{Lemma}
\aliascntresetthe{lemm}
\theoremstyle{definition}
\newtheorem{defn}[thm]{Definition}
\theoremstyle{remark}
\newtheorem{remark}[thm]{Remark}
\allowdisplaybreaks

\title{The Projective Reed--Muller Code from Quartic Curves in the Projective Plane}
\author{Xu Zhuang\thanks{Department of Mathematics, University of California,
Irvine, 340 Rowland Hall, Irvine, CA 92697-3875, USA.
Email: \href{mailto:xzhuang8@uci.edu}{\texttt{xzhuang8@uci.edu}}.}}
\date{October 14, 2026}
\hypersetup{
  pdftitle={The Projective Reed--Muller Code from Quartic Curves in the Projective Plane},
  pdfauthor={Xu Zhuang}
}

\begin{document}

\maketitle
\begin{abstract}
We study the Hamming weight enumerator of the projective Reed--Muller
code \(\cC_{2,4}\) over \(\F_q\), with \(\operatorname{char}(\F_q)>3\).
We determine the geometrically reducible and singular contributions by
combining incidence counts with counts of plane linear series on a fixed
normalization.  The line and conic configurations and the genus zero
cases give explicit polynomials.  The line--cubic and genus one cases
give formulas involving class numbers.  The genus two contribution is a finite sum
over hyperelliptic models, using their point counts over \(\F_q\) and
\(\F_{q^2}\) and their rational Weierstrass points.  The remaining smooth
contribution is the weighted point count distribution of
nonhyperelliptic genus three curves.  
\end{abstract}
\section{Introduction}
The Hamming weight enumerator of a linear code \(C\) of length \(\ell\)
over \(\F_q\) is
\[
 W_C(X,Y)=\sum_{c\in C}X^{\ell-\wt(c)}Y^{\wt(c)},
\]
where \(\wt(c)\) is the number of nonzero coordinates of \(c\).
The projective Reed--Muller code \(\cC_{m,d}\) is obtained by evaluating
the homogeneous degree \(d\) forms in \(m+1\) variables at chosen vector
representatives of the points of \(\PP^m(\F_q)\).  Changing these
representatives gives a monomially equivalent code and preserves its
weight enumerator
\cite{lachaud1990parameters,sorensen1991projective,kaplan2018weight}.
The zero coordinates of the codeword defined by a form \(F\) are exactly
the points of \(V(F)(\F_q)\).  In the injective range, computing the
weight enumerator therefore amounts to counting equations with each
possible number of rational zeros.

We consider ternary quartic forms.  Throughout the paper, let \(q=p^r\)
with \(p>3\), fix \(\overline{\F}_q\), and write \(\Frob_q\) for the
\(q\) power Frobenius.  Since \(q\ge5\), projective evaluation in degree
four is injective.  Thus \(\cC_{2,4}\) has length \(q^2+q+1\), dimension
15, and weight enumerator
\[
 W_{\cC_{2,4}}(X,Y)
 =\sum_{F\in\F_q[x_0,x_1,x_2]_4}
 X^{\#V(F)(\F_q)}Y^{q^2+q+1-\#V(F)(\F_q)}.
\]
Here \(V(0)=\PP^2\).  Nonzero scalar multiples of a form are distinct
codewords, although they define the same plane curve.

We call a form reducible if it factors over \(\overline{\F}_q\).
For an absolutely irreducible singular quartic, the normalization has
genus two, one, or zero.  A smooth plane quartic is a nonhyperelliptic
genus three curve in its canonical embedding.  Let
\(\mathcal M_3^{\nhyp}(\F_q)\) denote representatives of the
\(\F_q\)-isomorphism classes of such smooth curves, and set
\[
 S_q(n)=\sum_{\substack{[C]\in\mathcal M_3^{\nhyp}(\F_q)\\
                         \#C(\F_q)=n}}
              \frac1{\#\Aut_{\F_q}(C)}.
\]
The Hasse--Weil bound gives \(S_q(n)=0\) unless
\(\lvert n-(q+1)\rvert\le6\sqrt q\).

\begin{thm}\label{thm:intro-structure}
The weight enumerator has the decomposition
\[
\begin{aligned}
 W_{\cC_{2,4}}(X,Y)
 ={}&X^{q^2+q+1}+W_{\mathrm{red}}(X,Y)+W_{\mathrm{irr,sing}}(X,Y)\\
 &+(q-1)\#\PGL_3(\F_q)
      \sum_{n=0}^{q^2+q+1}S_q(n)X^nY^{q^2+q+1-n}.
\end{aligned}
\]
The reducible and absolutely irreducible singular contributions are
determined by explicit polynomials, formulas involving class numbers,
and finite sums over hyperelliptic models of genus two curves, as
summarized below.  The smooth contribution is expressed in terms of
\(S_q(n)\), which remains unevaluated for general \(q\).
\end{thm}

\begin{center}
\small
\renewcommand{\arraystretch}{1.15}
\begin{tabular}{@{}T{0.25\textwidth}T{0.39\textwidth}T{0.29\textwidth}@{}}
\toprule
Contribution & Arithmetic data & Form of the answer\\
\midrule
Line and conic configurations & Finite plane incidence data & Explicit polynomials\\[2pt]
Line--cubic and genus one cases & Elliptic traces and rational two-torsion & Formulas involving class numbers\\[2pt]
Genus zero normalization & Frobenius action on singular fibers & Explicit polynomials\\[2pt]
Genus two normalization & Point counts over \(\F_q,\F_{q^2}\) and rational Weierstrass points & Finite hyperelliptic model sums\\[2pt]
Smooth quartics & Nonhyperelliptic genus three point distribution & Weighted distribution \(S_q(n)\)\\
\bottomrule
\end{tabular}
\end{center}
Theorem~\ref{thm:full-weight-enumerator} gives the expanded formula and
identifies each coefficient.  The factor \((q-1)\#\PGL_3(\F_q)\) in the
smooth term accounts for scalar equations and projective coordinates,
with automorphisms accounted for in \(S_q(n)\).

For absolutely irreducible singular quartics, the common method is to fix
the normalization \(C\) and count base point free degree four plane linear
series that give birational maps onto their images.
For each singularity profile, Proposition~\ref{prop:irr-strata-contrib}
converts this count into a count of equations by multiplying by
\((q-1)\#\PGL_3(\F_q)/\#\Aut_{\F_q}(C)\).
The Frobenius action on the singular fibers determines the change in
rational point count from \(C\) to its plane model.  Thus the geometric
count and its contribution to the enumerator are treated in the same way
in each genus.

Two parameter spaces make this method particularly concrete.
For a genus two curve \(C\), let \(\operatorname{Weier}(C)(\F_q)\)
denote its rational Weierstrass points, the rational ramification points
of the hyperelliptic map.  The plane series on \(C\) are parametrized by
\((\Sym^2 C)(\F_q)\setminus|K_C|(\F_q)\).
The excluded canonical divisors give double covers of conics, while the
remaining divisors are the singular fibers of quartic models.
Proposition~\ref{prop:g2-fixed-counts} expresses the split node, cusp,
and nonsplit node counts using only \(\#C(\F_q)\),
\(\#C(\F_{q^2})\), and \(\#\operatorname{Weier}(C)(\F_q)\).
For genus one, projections of an elliptic normal quartic give the
parameter space.  Proposition~\ref{prop:g1-projection-geometry}
relates the projection centers to degree two divisors and divisor pairs,
identifies the exceptional conic maps, and determines the singular
fibers.  Counting these data gives the trace shift formulas of
Proposition~\ref{prop:g1-fixed-counts}.

The divisor constructions already occur in Bergstr\"om's classification
and count of singular quartics for a sieve used to study the cohomology
of genus three moduli spaces.  His genus one criterion identifies the
degree two maps by the linear equivalence of the two degree two divisors
\cite[Lemma~13.7]{bergstrom2008cohomology}.
We use his genus zero counts for two and three geometric singular points
\cite[Sections~12.3 and~13.4]{bergstrom2008cohomology}.
For one geometric singular point we give direct derivations using linear
series and intrinsic branch data.  We also give direct incidence and
projection counts for the reducible and genus one cases.
We organize these constructions and counts by fixed normalization and
trace shift to obtain the coefficients required by this weight enumerator.
The line--cubic contribution and the arithmetic evaluation of the elliptic
sums use the formulas of Kaplan and Kaplan--Petrow
\cite{kaplan2018weight,kaplanpetrow2017trace}.

This continues earlier work on exact weight enumerators of projective
Reed--Muller codes.  Elkies computed the weight enumerators for quadrics
in projective spaces of arbitrary dimension and for cubics in
\(\PP^2\) and \(\PP^3\) \cite{elkies2006linear}.
Aubry obtained related results for projective quadrics
\cite{aubry1992reed}.
Kaplan studied the weight enumerators of the projective and affine
Reed--Muller codes obtained from plane cubics \cite{kaplan2018weight}.
Kaplan--Matei used the second weight enumerator to count pairs of plane
cubics with prescribed rational intersections \cite{kaplan2021counting}.

The smooth contribution also connects this problem with the study of
rational points on genus three curves.  Lauter proved that over every
finite field there is a genus three curve whose number of rational points
lies within three of the Serre--Weil upper or lower bound
\cite{lauter2002maximum}.  The alternative between the two bounds is
related to Serre's obstruction.  A principally polarized abelian threefold
that is geometrically a Jacobian may require a quadratic twist to become
a Jacobian over \(\F_q\), and this twist changes the sign of its
Frobenius trace.  Ritzenthaler computed this obstruction explicitly and
applied it to the existence of optimal genus three curves
\cite{ritzenthaler2010serre}.

This obstruction also appears in the distribution of rational point counts.
Quadratic twisting pairs elliptic curves, and likewise hyperelliptic
curves, with opposite Frobenius traces.  For smooth plane quartics,
the quadratic twist of the Jacobian need not be a Jacobian over \(\F_q\),
so there is no analogous pairing.  Lercier, Ritzenthaler, Rovetta,
Sijsling, and Smith gave a heuristic explanation of the observed
asymmetry in the trace distribution and its relation to Serre's
obstruction \cite{lercier2019distributions}.
More recently, Bergstr\"om, Howe, Lorenzo Garc\'ia, and Ritzenthaler
refined Katz--Sarnak theory by determining limits of moments that measure
this asymmetry for nonhyperelliptic curves of genus at least three
\cite{bergstrom2025refinements}.
Our formulas isolate the distribution \(S_q(n)\) from the reducible and
singular contributions, making explicit which part of the quartic
enumerator depends on these questions about smooth curves.

Section~2 counts the reducible quartics by their geometric components
and Frobenius action.  The incidence tables are collected in the appendix.
Section~3 applies fixed normalization counting in genus two, one, and zero.
Section~4 assembles the enumerator, evaluates the elliptic sums by class
numbers, and writes the genus two contribution as a finite sum over all
squarefree polynomials of degree five or six.
SageMath code and saved outputs are maintained in a companion
repository.\footnote{\url{https://github.com/zxmath/quartic-weight-enumerator-calculations}.}
\section{Reducible quartic forms}
\label{sec:reducible}

We count nonzero quartic forms that factor over \(\overline{\F}_q\).
The point count depends only on the reduced union of the components.
For the line and conic configurations, we count equations up to
multiplication by a nonzero scalar using incidence arguments.
We treat the line and cubic case using the affine cubic weight enumerator.
The full incidence tables are in
\Cref{app:reducible-tables}.  We describe the counting methods below.

\subsection{Line and conic configurations}

For quartics with four line factors, the configurations and their counts
are listed in \Cref{tab:fourline}.  We check the total independently by
choosing Frobenius orbits of lines, allowing repetitions.
Let \(n_d\) be the number of such orbits of length \(d\).
There are \(q^{2d}+q^d+1\) lines over \(\F_{q^d}\), since lines are
points of the dual plane.  Removing the lines defined over proper
subfields and dividing by \(d\) gives
\(n_1=q^2+q+1\), \(n_2=(q^4-q)/2\),
\(n_3=(q^6+q^3-q^2-q)/3\), and \(n_4=(q^8-q^2)/4\).
The orbit lengths, counted with repetition, form one of the partitions
\(1+1+1+1\), \(2+1+1\), \(2+2\), \(3+1\), or \(4\).
These give the five terms in the total
\[
\begin{aligned}
&\binom{n_1+3}{4}+n_2\binom{n_1+1}{2}
 +\binom{n_2+1}{2}+n_1n_3+n_4\\
&=q^8+q^7+2q^6+2q^5+3q^4+2q^3+2q^2+q+1.
\end{aligned}
\]
This agrees with the sum of all 21 rows of \Cref{tab:fourline}.

There are \(q^5-q^2\) smooth conics over \(\F_q\), each with \(q+1\)
rational points.  A fixed smooth conic has \(q+1\) rational tangent
lines, \(q(q+1)/2\) rational lines meeting it in two rational points,
and \(q(q-1)/2\) rational lines meeting it in no rational point.
These counts give \Cref{tab:double-line-conic} directly.
For two distinct rational lines, choose an unordered pair of the
required types and separate the pairs whose intersection lies on the
conic.  Inclusion and exclusion then gives \Cref{tab:lines-conic}.

A pair of nonrational conjugate lines has just one rational point,
their intersection.  Each rational point lies on \((q^2-q)/2\) such
pairs.  For a fixed smooth conic, choosing this point on or off the
conic gives \Cref{tab:glc}.

For two distinct smooth conics, we need only the number \(r\) of
\emph{distinct} rational intersection points.  B\'ezout gives
\(0\le r\le4\).  A tangency contributes one point to \(r\).
If both conics are defined over \(\F_q\), their union has
\(2q+2-r\) rational points.  If they are conjugate over \(\F_{q^2}\),
a rational point on either conic also lies on its conjugate, so their
union has exactly \(r\) rational points.

We count the two families by the same procedure: choose some rational
intersection points, count pairs through them, and remove pairs with
more intersection points.  For \(0\leq i\leq4\), let \(a_i\) be the number of subsets
of \(i\) points in \(\PP^2(\F_q)\) with no three collinear, and let
\(c_i(q)\) be the number of smooth conics through any one such subset.
The elementary counts are
\begingroup
\setlength{\LTpre}{6pt}
\setlength{\LTpost}{6pt}
\renewcommand{\LTcaptype}{unnumtable}
\begin{longtable}{@{}ccc@{}}
\toprule
\(i\) & \(a_i\) & \(c_i(q)\)\\
\midrule
\endfirsthead
\toprule
\(i\) & \(a_i\) & \(c_i(q)\)\\
\midrule
\endhead
\midrule
\multicolumn{3}{r}{Continued on next page}\\
\endfoot
\bottomrule
\endlastfoot
0 & \(1\) & \(q^5-q^2\)\\[3pt]
1 & \(q^2+q+1\) & \(q^4-q^2\)\\[3pt]
2 & \(\displaystyle\binom{q^2+q+1}{2}\) & \(q^3-q^2\)\\[7pt]
3 & \(\displaystyle\frac{(q^2+q+1)(q^2+q)q^2}{6}\) & \((q-1)^2\)\\[7pt]
4 & \(\displaystyle\frac{(q^2+q+1)(q^2+q)q^2(q-1)^2}{24}\) & \(q-2\).\\
\end{longtable}
\endgroup
\noindent
The values of \(a_i\) follow by choosing the points successively.
Since \(\PGL_3(\F_q)\) acts transitively on these subsets, counting
incidences with smooth conics gives
\(a_i c_i(q)=(q^5-q^2)\binom{q+1}{i}\).

\begin{figure}[!htbp]
\centering
\begingroup
\pgfmathsetmacro{\conicIntersectionX}{sqrt((1/0.85^2-1/1.3^2)/
  (1/(1.15^2*0.85^2)-1/(1.8^2*1.3^2)))}
\pgfmathsetmacro{\conicIntersectionY}{0.85*sqrt(1-\conicIntersectionX^2/1.8^2)}
\begin{tikzpicture}[
  x=1cm,y=1cm,font=\small,
  firstconic/.style={draw=black,line width=1.05pt},
  conjugateconic/.style={draw=black,line width=0.8pt,dash pattern=on 4pt off 2.8pt},
  rationalpoint/.style={circle,fill=black,inner sep=1.5pt},
  conjugatepoint/.style={circle,draw=black,fill=white,line width=0.7pt,inner sep=1.65pt},
  every node/.style={outer sep=2pt}]
\path[use as bounding box] (-3.7,-0.95) rectangle (11.2,4.15);
\foreach \panel/\shift in {1/0,2/7.5}{
\begin{scope}[xshift=\shift cm]
  \ifnum\panel=1
    \node at (0,3.75) {Two rational intersection points};
  \else
    \node at (0,3.75) {Four rational intersection points};
  \fi
  \draw[firstconic] (0,1.45) ellipse[x radius=1.8,y radius=0.85];
  \draw[conjugateconic] (0,1.45) ellipse[x radius=1.15,y radius=1.3];
  \node[anchor=west] at (1.91,1.45) {\(Q\)};
  \node at (0,3.03) {\(\Frob_q(Q)\)};
  \ifnum\panel=1
    \draw[<->,line width=0.6pt,shorten <=4pt,shorten >=4pt]
      (-\conicIntersectionX,{1.45+\conicIntersectionY})
      .. controls (-3,3) and (-3,-0.1) ..
      (-\conicIntersectionX,{1.45-\conicIntersectionY});
    \node[anchor=east] at (-2.48,1.45) {\(\Frob_q\)};
    \foreach \sign in {-1,1}{
      \node[conjugatepoint] at
        (-\conicIntersectionX,{1.45+\sign*\conicIntersectionY}) {};
      \node[rationalpoint] at
        (\conicIntersectionX,{1.45+\sign*\conicIntersectionY}) {};
    }
    \node at (0,-0.3) {\(\#(Q\cup\Frob_q(Q))(\F_q)=2\)};
    \node at (0,-0.8) {Choices of two rational points: \(\binom{2}{2}=1\)};
  \else
    \foreach \side in {-1,1}{
      \foreach \sign in {-1,1}{
        \node[rationalpoint] at
          ({\side*\conicIntersectionX},{1.45+\sign*\conicIntersectionY}) {};
      }
    }
    \node at (0,-0.3) {\(\#(Q\cup\Frob_q(Q))(\F_q)=4\)};
    \node at (0,-0.8) {Choices of two rational points: \(\binom{4}{2}=6\)};
  \fi
\end{scope}
}
\end{tikzpicture}
\endgroup
\caption{Two schematic pairs of conjugate smooth conics
\(Q,\Frob_q(Q)\), with \(Q\) defined over \(\F_{q^2}\) but not
over \(\F_q\), drawn over \(\overline{\F}_q\).
Frobenius exchanges the two conics in each pair.
Filled intersection points are rational, and the
two open points are exchanged by Frobenius.  The unions therefore
have two and four rational points, respectively.  Choosing two rational
intersection points counts the left pair once and the right pair six
times.}
\label{fig:reducible-incidences}
\end{figure}
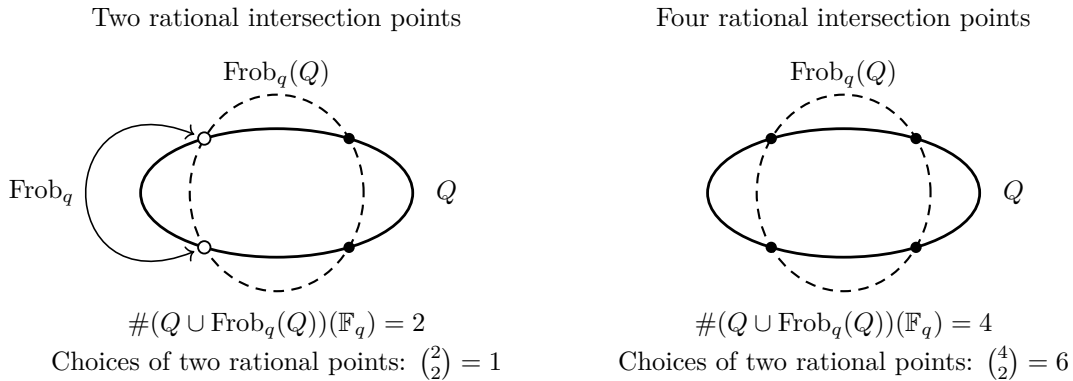

Treat each family separately, and let \(N_r\) count its unordered
pairs with exactly \(r\) rational intersection points.  A pair counted by
\(N_r\) contains \(\binom ri\) subsets of \(i\) such points.
Alternatively, choose the subset first, in \(a_i\) ways.
For two rational conics, choose two of the \(c_i(q)\) smooth conics
through it.  For a conjugate pair, there are \(c_i(q^2)-c_i(q)\)
nonrational smooth conics through it.  Dividing by two identifies each
conic with its conjugate.  Thus, for \(i=0,\ldots,4\),
\[
 \sum_{r=i}^{4}\binom ri N_r
 =\begin{cases}
   a_i\binom{c_i(q)}2,&\text{two distinct rational conics},\\[3pt]
   \dfrac{a_i}{2}\bigl(c_i(q^2)-c_i(q)\bigr),
       &\text{two conjugate conics}.
  \end{cases}
\]
These equations determine the counts from \(r=4\) down to \(r=0\).
The equation for \(i=4\) gives \(N_4\) directly.  For \(i=3\), the
left side is \(N_3+4N_4\): a pair with four rational intersection
points contains four subsets of three.  For example, in the conjugate
case this gives \(N_4=a_4(q^2-q)/2\) and
\(N_3=a_3\bigl((q^2-1)^2-(q-1)^2\bigr)/2-4N_4\).
Next, the equation for \(i=2\) gives \(N_2\) after subtracting
\(3N_3+6N_4\), and the cases \(i=1,0\) follow in the same way.
This gives the distinct-conic rows in \Cref{tab:two-conics,tab:GC}.
The \(q^5-q^2\) double smooth conics each have \(q+1\) rational
points and are added separately to \Cref{tab:two-conics}.
Figure~\ref{fig:reducible-incidences} illustrates the repeated counting
when two intersection points are selected.

\long\def\FourLineCountingTable{%
\begingroup
\small
\renewcommand{\arraystretch}{1.1}
\setlength{\tabcolsep}{5pt}
\setlength{\LTpre}{6pt}
\setlength{\LTpost}{6pt}

\begin{longtable}{@{}
    T{6.6cm}
    Q{1.9cm}
    Q{6.4cm}
@{}}
\caption{Four line quartics over $\F_q$}\label{tab:fourline}\\
\toprule
\textbf{Type} & \textbf{Number of \mbox{$\F_q$-points}} & \textbf{Number of quartics} \\
\midrule
\endfirsthead

\multicolumn{3}{c}{\tablename\ \thetable\ (continued from previous page)} \\
\toprule
\textbf{Type} & \textbf{Number of \mbox{$\F_q$-points}} & \textbf{Number of quartics} \\
\midrule
\endhead

\midrule
\multicolumn{3}{r}{Continued on next page} \\
\endfoot

\bottomrule
\endlastfoot

Four identical lines over $\F_q$
& $q+1$
& $q^2+q+1$ \\

A triple line and one other line over $\F_q$
& $2q+1$
& $(q^2+q+1)(q^2+q)$ \\

Two distinct double lines over $\F_q$
& $2q+1$
& $\dfrac{(q^2+q+1)(q^2+q)}{2}$ \\

One pair of conjugate double lines over $\F_{q^2}\setminus \F_q$
& $1$
& $\dfrac{(q^2+q+1)(q^2-q)}{2}$ \\

Two distinct lines over $\F_q$ intersecting a double line over $\F_q$ at different points
& $3q$
& $(q^2+q+1)\binom{q+1}{2}q^2$ \\

Two distinct lines over $\F_q$ intersecting a double line over $\F_q$ at the same point
& $3q+1$
& $3(q^2+q+1)\binom{q+1}{3}$ \\

One rational double line intersecting two Galois conjugate lines at different points
& $q+2$
& $\dfrac{(q^2-q)q^2}{2}(q^2+q+1)$ \\

One rational double line intersecting two Galois conjugate lines at the same point
& $q+1$
& $\dfrac{(q^2-q)(q+1)(q^2+q+1)}{2}$ \\

Four distinct rational lines intersecting at the same point
& $4q+1$
& $(q^2+q+1)\binom{q+1}{4}$ \\

Three rational lines intersecting at the same point with one additional line
& $4q-1$
& $(q^2+q+1)\binom{q+1}{3}q^2$ \\

Four distinct rational lines intersecting at six points
& $4q-2$
& $\dfrac{(q^2+q+1)(q^2+q)q^2(q^2-2q+1)}{24}$ \\

Two rational lines and two Galois conjugate lines intersecting at the same point
& $2q+1$
& $(q^2+q+1)\dfrac{q^2-q}{2}\binom{q+1}{2}$ \\

One rational line passing through the intersection point of two Galois conjugate lines and one rational line not passing through the intersection point
& $2q+1$
& $(q^2+q+1)\dfrac{q^2-q}{2}(q+1)q^2$ \\

Two Galois conjugate lines and two rational lines intersecting at different points
& $2q+2$
& $(q^2+q+1)\dfrac{q^2-q}{2}\binom{q^2}{2}$ \\

One rational line with three $\Gal(\F_{q^3}/\F_q)$ conjugate lines through a rational point not on the rational line
& $q+2$
& $\dfrac{q^3-q}{3}(q^2+q+1)q^2$ \\

One rational line with three $\Gal(\F_{q^3}/\F_q)$ conjugate lines through a rational point on the rational line
& $q+1$
& $(q^2+q+1)\dfrac{q^3-q}{3}(q+1)$ \\

One rational line with three conjugate lines that do not all pass through one point
& $q+1$
& $(q-1)\dfrac{q^5-q^3}{3}(q^2+q+1)$ \\

Two pairs of $\Gal(\F_{q^2}/\F_q)$ conjugate lines intersecting at the same point
& $1$
& $(q^2+q+1)\binom{\frac{q^2-q}{2}}{2}$ \\

Two pairs of $\Gal(\F_{q^2}/\F_q)$ conjugate lines intersecting at two points
& $2$
& $\binom{q^2+q+1}{2}\left(\dfrac{q^2-q}{2}\right)^2$ \\

Four $\Gal(\F_{q^4}/\F_q)$ conjugate lines
& $1$
& $(q^2+q+1)\dfrac{q^4-q^2}{4}$ \\

Four Galois conjugate lines that do not all pass through one point
& $0$
& $\dfrac{(q^8-q^2)-(q^2+q+1)(q^4-q^2)}{4}$ \\

\end{longtable}
\normalsize
\endgroup
}

\long\def\DoubleLineConicCountingTable{%
\begingroup
\small
\renewcommand{\arraystretch}{1.1}
\setlength{\tabcolsep}{5pt}
\setlength{\LTpre}{6pt}
\setlength{\LTpost}{6pt}

\begin{longtable}{@{}
    T{6.6cm}
    Q{1.9cm}
    Q{6.4cm}
@{}}
\caption{A double line and a conic}\label{tab:double-line-conic}\\
\toprule
\textbf{Type} & \textbf{Number of \mbox{$\F_q$-points}} & \textbf{Number of quartics} \\
\midrule
\endfirsthead

\multicolumn{3}{c}{\tablename\ \thetable\ (continued from previous page)} \\
\toprule
\textbf{Type} & \textbf{Number of \mbox{$\F_q$-points}} & \textbf{Number of quartics} \\
\midrule
\endhead

\midrule
\multicolumn{3}{r}{Continued on next page} \\
\endfoot

\bottomrule
\endlastfoot

\small
The double line is tangent to the conic
& $2q+1$
& $(q+1)(q^5-q^2)$ \\

The double line passes through two rational points on the conic
& $2q$
& $\dfrac{(q+1)q(q^5-q^2)}{2}$ \\

The double line passes through two Galois conjugate points on the conic
& $2q+2$
& $\dfrac{(q-1)q(q^5-q^2)}{2}$ \\

\end{longtable}
\normalsize
\endgroup
}

\long\def\TwoLinesConicCountingTable{%
\begingroup
\small
\renewcommand{\arraystretch}{1.1}
\setlength{\tabcolsep}{5pt}
\setlength{\LTpre}{6pt}
\setlength{\LTpost}{6pt}

\begin{longtable}{@{}
    T{6.6cm}
    Q{1.9cm}
    Q{6.4cm}
@{}}
\caption{A pair of distinct rational lines and a smooth conic}\label{tab:lines-conic}\\
\toprule
\textbf{Type} & \textbf{Number of \mbox{$\F_q$-points}} & \textbf{Number of quartics} \\
\midrule
\endfirsthead

\multicolumn{3}{c}{\tablename\ \thetable\ (continued from previous page)} \\
\toprule
\textbf{Type} & \textbf{Number of \mbox{$\F_q$-points}} & \textbf{Number of quartics} \\
\midrule
\endhead

\midrule
\multicolumn{3}{r}{Continued on next page} \\
\endfoot

\bottomrule
\endlastfoot

\small
Both lines tangent to the conic
& $3q$
& $\dfrac{(q+1)q(q^5-q^2)}{2}$ \\

One line tangent to the conic, the other line passing through the tangent point
& $3q$
& $(q+1)q(q^5-q^2)$ \\

One line tangent to the conic, the other line passing through two rational points on the conic (not the tangent point)
& $3q-1$
& $\dfrac{(q+1)q(q-1)(q^5-q^2)}{2}$ \\

One line is tangent to the conic, and the other passes through two Galois conjugate points on the conic
& $3q+1$
& $\dfrac{(q-1)q(q+1)(q^5-q^2)}{2}$ \\

Both lines pass through two rational points on the conic, and their intersection point lies on the conic
& $3q-1$
& $\dfrac{(q+1)q(q-1)(q^5-q^2)}{2}$ \\

Both lines pass through two rational points on the conic, and their intersection point does not lie on the conic
& $3q-2$
& $\dfrac{(q+1)q(q-1)(q-2)(q^5-q^2)}{8}$ \\

One line passes through two rational points on the conic, the other line passes through two Galois conjugate points on the conic
& $3q$
& $\dfrac{(q-1)q^2(q+1)(q^5-q^2)}{4}$ \\

Both lines pass through two Galois conjugate points on the conic
& $3q+2$
& $\dfrac{(q+1)q(q-1)(q-2)(q^5-q^2)}{8}$ \\

\end{longtable}
\normalsize
\endgroup
}

\long\def\TwoConicsCountingTable{%
\begingroup
\small
\renewcommand{\arraystretch}{1.1}
\setlength{\tabcolsep}{5pt}
\setlength{\LTpre}{6pt}
\setlength{\LTpost}{6pt}

\begin{longtable}{@{}
    T{6.3cm}
    Q{1.9cm}
    Q{6.7cm}
@{}}
\caption{Two smooth conics}\label{tab:two-conics}\\
\toprule
\textbf{Type} & \textbf{Number of \mbox{$\F_q$-points}} & \textbf{Number of quartics} \\
\midrule
\endfirsthead

\multicolumn{3}{c}{\tablename\ \thetable\ (continued from previous page)} \\
\toprule
\textbf{Type} & \textbf{Number of \mbox{$\F_q$-points}} & \textbf{Number of quartics} \\
\midrule
\endhead

\midrule
\multicolumn{3}{r}{Continued on next page} \\
\endfoot

\bottomrule
\endlastfoot

\small
A double smooth conic
& $q+1$
& $q^5-q^2$ \\

Two distinct smooth conics intersecting at exactly four rational points
& $2q-2$
& $\dfrac{(q+1)(q-1)(q-2)(q-3)q(q^5-q^2)}{48}$ \\

Two distinct smooth conics intersecting at exactly three rational points
& $2q-1$
& $\dfrac{(q+1)(q-1)(q-2)q(q^5-q^2)}{4}$ \\

Two distinct smooth conics intersecting at exactly two rational points
& $2q$
& $\dfrac{(q^3-2q^2+7q-8)(q+1)q(q^5-q^2)}{8}$ \\

Two distinct smooth conics intersecting at exactly one rational point
& $2q+1$
&  $\dfrac{(2q^3+5q^2-6q+6)(q+1)(q-1)(q^5-q^2)}{12}$ \\

Two distinct smooth conics with no rational intersection points
& $2q+2$
& $\dfrac{(3q^2-q+2)(q-1)^2q(q^5-q^2)}{16}$ \\

\end{longtable}
\normalsize
\endgroup
}

\long\def\ConjugateLinesConicCountingTable{%
\begingroup
\small
\renewcommand{\arraystretch}{1.1}
\setlength{\tabcolsep}{5pt}
\setlength{\LTpre}{6pt}
\setlength{\LTpost}{6pt}

\begin{longtable}{@{}
    T{6.6cm}
    Q{1.9cm}
    Q{6.4cm}
@{}}
\caption{A pair of Galois conjugate lines with a smooth conic}\label{tab:glc}\\
\toprule
\textbf{Type} & \textbf{Number of \mbox{$\F_q$-points}} & \textbf{Number of quartics} \\
\midrule
\endfirsthead

\multicolumn{3}{c}{\tablename\ \thetable\ (continued from previous page)} \\
\toprule
\textbf{Type} & \textbf{Number of \mbox{$\F_q$-points}} & \textbf{Number of quartics} \\
\midrule
\endhead

\midrule
\multicolumn{3}{r}{Continued on next page} \\
\endfoot

\bottomrule
\endlastfoot

\small
The intersection point of the conjugate lines lies on the conic
& $q+1$
& $\dfrac{q(q-1)(q+1)(q^5-q^2)}{2}$ \\

The intersection point of the conjugate lines does not lie on the conic
& $q+2$
& $\dfrac{(q^2-q)q^2(q^5-q^2)}{2}$ \\

\end{longtable}
\normalsize
\endgroup
}

\long\def\ConjugateConicsCountingTable{%
\begingroup
\small
\renewcommand{\arraystretch}{1.1}
\setlength{\tabcolsep}{5pt}
\setlength{\LTpre}{6pt}
\setlength{\LTpost}{6pt}
\begin{longtable}{@{}
    T{6.6cm}
    Q{1.9cm}
    Q{6.4cm}
@{}}
\caption{Two Galois conjugate conics}\label{tab:GC}\\
\toprule
\textbf{Type} & \textbf{Number of \mbox{$\F_q$-points}} & \textbf{Number of quartics} \\
\midrule
\endfirsthead

\multicolumn{3}{c}{\tablename\ \thetable\ (continued from previous page)} \\
\toprule
\textbf{Type} & \textbf{Number of \mbox{$\F_q$-points}} & \textbf{Number of quartics} \\
\midrule
\endhead

\midrule
\multicolumn{3}{r}{Continued on next page} \\
\endfoot

\bottomrule
\endlastfoot

\small
Two Galois conjugate conics intersecting at exactly one rational point
& $1$
& \makecell[l]{$\dfrac{1}{12}(2q^3+q^2+6)(q^2+q+1)$\\
$(q+1)(q-1)^2q^2$} \\

Two Galois conjugate conics intersecting at exactly two rational points
& $2$
& \makecell[l]{$\dfrac{1}{8}(q^3-2q^2+5q-2)(q^2+q+1)$\\
$(q+1)(q-1)q^3$} \\

Two Galois conjugate conics intersecting at exactly three rational points
& $3$
& $\begin{aligned}
\frac{1}{4}(q^2+q+1)(q+1)(q-1)^2q^4
\end{aligned}$ \\

Two Galois conjugate conics intersecting at exactly four rational points
& $4$
& $\begin{aligned}
\frac{1}{48}(q^2+q+1)(q+1)(q-1)^3q^4
\end{aligned}$ \\

Two Galois conjugate conics with no rational intersection points
& $0$
& $\begin{aligned}
\frac{1}{16}(3q^3-3q-4)(q^2+q+1)(q-1)^2q^3
\end{aligned}$ \\

\end{longtable}
\normalsize
\endgroup
}

\subsection{A line and an absolutely irreducible cubic}

Suppose that a quartic form has geometric factorization \(F=LG\), where
\(L\) is linear and \(G\) is absolutely irreducible of degree three.  The
components have different degrees, so each is fixed by Frobenius.  Thus both
components are defined over \(\F_q\), and the line component is unique.
Fix an equation for this line and choose coordinates in which it is
\(X_0=0\).  Setting \(X_0=1\) gives a bijection between the homogeneous
cubic factors and affine polynomials \(g(x,y)=G(1,x,y)\) for which \(G\)
is absolutely irreducible.
These polynomials give the corresponding part of the affine Reed--Muller code
\(\cC^A_{2,3}\subseteq\F_q^{q^2}\).
The rational zero set of \(X_0G\) is the disjoint union of the \(q+1\)
points on \(X_0=0\) and the affine zero set of \(g\).  Thus, for a fixed
line, the quartic weight equals the affine weight of \(g\).
The affine enumerators below count all nonzero scalar equations for each
cubic type.

We use the following two results on affine cubics.

\begin{prop}[{\cite[Section~5.1, Proposition~4]{kaplan2018weight}}]
The contributions from cuspidal, split nodal, and nonsplit nodal absolutely
irreducible affine cubics are, respectively,
\begin{align*}
W_3^{A,\mathrm{cusp}}(X,Y)
={}&(q-1)(q^3-q)q^2
\bigg(
 \frac{(q+1)(q-1)}3X^{q+1}Y^{q^2-q-1}
 +\frac{q^2-q+4}{2}X^qY^{q^2-q}\\
&\hspace{3.2cm}
 +(2q-1)X^{q-1}Y^{q^2-q+1}
 +\frac{(q-1)(q-2)}6X^{q-2}Y^{q^2-q+2}
\bigg),\\
W_3^{A,\mathrm{split}}(X,Y)
={}&\frac{(q-1)(q^3-q)(q^3-q^2)}2
\bigg(
 \frac{q(q+1)}3X^qY^{q^2-q}
 +\frac{q^2-q+6}{2}X^{q-1}Y^{q^2-q+1}\\
&\hspace{3.2cm}
 +(2q-3)X^{q-2}Y^{q^2-q+2}
 +\frac{(q-2)(q-3)}6X^{q-3}Y^{q^2-q+3}
\bigg),\\
W_3^{A,\mathrm{nonsplit}}(X,Y)
={}&\frac{(q-1)(q^3-q)(q^3-q^2)}2
\bigg(
 \frac{q(q-1)}3X^{q+2}Y^{q^2-q-2}
 +\frac{q(q-1)}2X^{q+1}Y^{q^2-q-1}\\
&\hspace{3.2cm}
 +(2q+1)X^qY^{q^2-q}
 +\frac{q(q-1)}6X^{q-1}Y^{q^2-q+1}
\bigg).
\end{align*}
\end{prop}

The four terms in each formula correspond, in order, to a line at infinity
meeting the projective cubic in \(i=0,1,2,3\) distinct rational points.  The
calculation of these line sections is given in
\cite[Section~5.1, Proposition~3 and equations~(4)--(7)]{kaplan2018weight}.

\begin{prop}[{\cite[Section~6,
Theorem~4]{kaplan2018weight}}]
For a pointed elliptic curve \((E,O)/\F_q\), put
\(t_E=q+1-\#E(\F_q)\) and define
\begin{align*}
 W^\alpha(X,Y,t)
={}&\frac{q^2+qt+t^2-q+t}{3}
       X^{q+1-t}Y^{q^2-q-1+t}
 +\frac{q^2-t^2+q+t+2}{2}
       X^{q-t}Y^{q^2-q+t}\\
 &+(q-t)X^{q-1-t}Y^{q^2-q+1+t}
 +\frac{(q-t)(q-t-1)}6
       X^{q-2-t}Y^{q^2-q+2+t}.
\end{align*}
Then the contribution from smooth absolutely irreducible affine cubics is
\begin{equation}\label{eq:affine-cubic-smooth}
 W_3^{A,\mathrm{sm}}(X,Y)
 =q^3(q-1)^2(q^2-1)
  \sum_{[E,O]/\F_q}
  \frac{W^\alpha(X,Y,t_E)}{\#\Aut_{\F_q}(E,O)}.
\end{equation}
\end{prop}

In \cite[Theorem~4]{kaplan2018weight}, the sum is grouped by trace and has
the additional condition \(t\not\equiv0\pmod3\).  This condition is a
typographical error: the last step of the proof sums over all elliptic curves, as
does Proposition~5 of that paper.  Formula~\eqref{eq:affine-cubic-smooth}
sums over all isomorphism classes and follows from Proposition~5 and the
proof of Theorem~4.

Since the line component of \(F\) is unique, the \(q^2+q+1\) choices of
line give disjoint sets of quartics.  The affine cubic enumerators already
count the nonzero scalar multiples of \(F\).  Thus quartics of this type
contribute
\begin{equation}\label{eq:line-cubic-contribution}
 (q^2+q+1)X^{q+1}
 \bigl(W_3^{A,\mathrm{cusp}}(X,Y)+W_3^{A,\mathrm{split}}(X,Y)
 +W_3^{A,\mathrm{nonsplit}}(X,Y)+W_3^{A,\mathrm{sm}}(X,Y)\bigr).
\end{equation}

\subsection{The reducible contribution}

The six line and conic factorization families in
Tables~\ref{tab:fourline}--\ref{tab:GC} are disjoint.
For each number \(z\) of rational zeros,
Table~\ref{tab:factorable-collected} gives their combined count \(R_z(q)\).
The incidence recurrences and the collection of these coefficients are
implemented in the accompanying code.
\par

\begingroup
\small
\begin{longtable}{@{}Q{0.13\linewidth}T{0.82\linewidth}@{}}
\caption{Collected coefficients for the six line and conic families.}
\label{tab:factorable-collected}\\
\toprule
Number of zeros $z$ & Projective coefficient $R_z(q)$\\
\midrule
\endfirsthead
\toprule
Number of zeros $z$ & Projective coefficient $R_z(q)$\\
\midrule
\endhead
\midrule
\multicolumn{2}{r}{\footnotesize Continued on the next page}\\
\endfoot
\bottomrule
\endlastfoot
$0$ & \adjustbox{max width=\linewidth}{$\displaystyle \frac{q^{4}(q-1)^{2}(q^{2}+q+1)(3q^{2}+1)}{16}$}\\[7pt]
$1$ & \adjustbox{max width=\linewidth}{$\displaystyle \frac{q(q-1)(q^{2}+q+1)(4q^{6}+2q^{5}-4q^{4}+10q^{3}+9q^{2}-9q+6)}{24}$}\\[7pt]
$2$ & \adjustbox{max width=\linewidth}{$\displaystyle \frac{q^{3}(q-1)(q+1)(q^{2}+q+1)(q^{3}-2q^{2}+6q-3)}{8}$}\\[7pt]
$3$ & \adjustbox{max width=\linewidth}{$\displaystyle \frac{q^{4}(q-1)^{2}(q+1)(q^{2}+q+1)}{4}$}\\[7pt]
$4$ & \adjustbox{max width=\linewidth}{$\displaystyle \frac{q^{4}(q-1)^{3}(q+1)(q^{2}+q+1)}{48}$}\\[7pt]
$q + 1$ & \adjustbox{max width=\linewidth}{$\displaystyle \frac{(q^{2}+q+1)(5q^{6}-5q^{5}-3q^{4}+16q^{3}-8q^{2}-5q+6)}{6}$}\\[7pt]
$q + 2$ & \adjustbox{max width=\linewidth}{$\displaystyle \frac{q^{3}(q-1)(q^{2}+q+1)(3q^{3}-3q^{2}+2q+5)}{6}$}\\[7pt]
$2 q - 2$ & \adjustbox{max width=\linewidth}{$\displaystyle \frac{q^{3}(q-3)(q-2)(q-1)^{2}(q+1)(q^{2}+q+1)}{48}$}\\[7pt]
$2 q - 1$ & \adjustbox{max width=\linewidth}{$\displaystyle \frac{q^{3}(q-2)(q-1)^{2}(q+1)(q^{2}+q+1)}{4}$}\\[7pt]
$2 q$ & \adjustbox{max width=\linewidth}{$\displaystyle \frac{q^{3}(q-1)(q+1)(q^{2}+q+1)(q^{3}-2q^{2}+7q-4)}{8}$}\\[7pt]
$2 q + 1$ & \adjustbox{max width=\linewidth}{$\displaystyle \frac{q(q+1)(q^{2}+q+1)(2q^{6}+q^{5}-14q^{4}+29q^{3}-9q^{2}-9q+18)}{12}$}\\[7pt]
$2 q + 2$ & \adjustbox{max width=\linewidth}{$\displaystyle \frac{q^{3}(q-1)^{2}(q^{2}+q+1)(3q^{3}-4q^{2}+7q+10)}{16}$}\\[7pt]
$3 q - 2$ & \adjustbox{max width=\linewidth}{$\displaystyle \frac{q^{3}(q-2)(q-1)^{2}(q+1)(q^{2}+q+1)}{8}$}\\[7pt]
$3 q - 1$ & \adjustbox{max width=\linewidth}{$\displaystyle q^{3}(q-1)^{2}(q+1)(q^{2}+q+1)$}\\[7pt]
$3 q$ & \adjustbox{max width=\linewidth}{$\displaystyle \frac{q^{3}(q+1)(q^{2}+q+1)(q^{3}-2q^{2}+7q-4)}{4}$}\\[7pt]
$3 q + 1$ & \adjustbox{max width=\linewidth}{$\displaystyle \frac{q(q-1)(q+1)(q^{2}+q+1)(q^{3}-q^{2}+1)}{2}$}\\[7pt]
$3 q + 2$ & \adjustbox{max width=\linewidth}{$\displaystyle \frac{q^{3}(q-2)(q-1)^{2}(q+1)(q^{2}+q+1)}{8}$}\\[7pt]
$4 q - 2$ & \adjustbox{max width=\linewidth}{$\displaystyle \frac{q^{3}(q-1)^{2}(q+1)(q^{2}+q+1)}{24}$}\\[7pt]
$4 q - 1$ & \adjustbox{max width=\linewidth}{$\displaystyle \frac{q^{3}(q-1)(q+1)(q^{2}+q+1)}{6}$}\\[7pt]
$4 q + 1$ & \adjustbox{max width=\linewidth}{$\displaystyle \frac{q(q-2)(q-1)(q+1)(q^{2}+q+1)}{24}$}\\[7pt]
\end{longtable}
\endgroup

Summing the seven disjoint factorization families gives
\begin{equation}
\begin{aligned}
 W_{\mathrm{red}}(X,Y)
 ={}&(q-1)\sum_{z}R_z(q)X^zY^{q^2+q+1-z}\\
 &+(q^2+q+1)X^{q+1}
   \bigl(W_3^{A,\mathrm{cusp}}(X,Y)+W_3^{A,\mathrm{split}}(X,Y)\\
 &\hspace{3.6cm}{}+W_3^{A,\mathrm{nonsplit}}(X,Y)
   +W_3^{A,\mathrm{sm}}(X,Y)\bigr),
\end{aligned}
\label{eq:factorable-enumerator}
\end{equation}
where \(z\) runs through the twenty entries of
Table~\ref{tab:factorable-collected}.  The first sum combines the
contributions of the six line and conic families, with a factor of
\(q-1\) for the nonzero scalar multiples.  The second is the line and cubic contribution
\eqref{eq:line-cubic-contribution}.
\par
\section{Absolutely irreducible quartics}

A smooth plane quartic is a nonhyperelliptic genus three curve embedded
by its complete canonical series, and every such curve has a canonical
plane model \cite[\S6.4]{eisenbud2024practice}.
Its contribution is the smooth term of Theorem~\ref{thm:intro-structure}.
We now count singular quartics by fixing their normalization, which has
genus two, one, or zero, and determining their plane linear series.

\subsection{Singularities, profiles, and series on a fixed normalization}
\label{subsec:profiles-fixed-normalization}

Let \(C_0\subset\PP^2_{\F_q}\) be an absolutely irreducible plane quartic,
and let \(\nu\colon C\longrightarrow C_0\) be its normalization.  Put
\(\overline C_0=C_0\otimes_{\F_q}\overline{\F}_q\).  We similarly write
\(\overline C=C\otimes_{\F_q}\overline{\F}_q\),
and let \(\overline\nu\colon\overline C\to\overline C_0\) be the base
change of \(\nu\).  It is again the normalization
\cite[Tag~0C3N]{stacksproject}.

Throughout this section, a \emph{branch} means a geometric branch, over
\(\overline{\F}_q\).  Intuitively, a branch at a point \(P\) is
one local irreducible piece of the curve passing through \(P\).

To describe this precisely, choose homogeneous coordinates \([X:Y:Z]\)
with \(P=[a:b:1]\), and let \(F(X,Y,Z)=0\) define \(\overline C_0\).
On the affine chart \(Z\ne0\), put \(u=X/Z-a\) and \(v=Y/Z-b\).
Then \(P=(0,0)\), and the curve has equation
\[
 f(u,v)=F(u+a,v+b,1)=0.
\]
We call this a \emph{local equation at \(P\)}: it is the equation of
the curve on this chart, with \(P\) moved to the origin.
The variables \(u,v\) are coordinates of the ambient plane.

To study branches, we regard \(f\) as an element of \(\overline{\F}_q[[u,v]]\),
the ring of formal power series.  Multiplying \(f\) by a unit, that is,
a power series with nonzero constant term, does not change the ideal
\((f)\).  Each irreducible factor \(f_B\), counted up to multiplication
by a unit, defines a branch \(B\) by the equation \(f_B=0\). See
\cite[Definition~16.3]{kunz2005plane}.
The normalization has one point above \(P\) for each branch
\cite[Theorem~16.14]{kunz2005plane}.
An absolutely irreducible curve can have several branches at a point,
because its defining polynomial can split in formal power series.

Define
\(\Sigma(C_0)=\operatorname{Sing}(\overline C_0)(\overline{\F}_q)\), and
let \(\mathcal B(C_0)=\overline\nu^{-1}(\Sigma(C_0))\).  Write
\(\pi_{C_0}\colon\mathcal B(C_0)\longrightarrow\Sigma(C_0)\) for the
restriction of \(\overline\nu\).  We identify each point of
\(\mathcal B(C_0)\) with the branch it represents.  Thus the number of
branches at \(P\) is \(r(P)=\#\pi_{C_0}^{-1}(P)\).
Both \(\Sigma(C_0)\) and \(\mathcal B(C_0)\) are finite sets of
geometric points, and \(\pi_{C_0}\) commutes with \(\Frob_q\).

Let \(M\) be a module over a ring \(R\).  Its \emph{length}, denoted by
\(\operatorname{length}_R(M)\), is the largest integer \(n\) for which
there is a chain of \(R\)-submodules
\[
 0=M_0\subsetneq M_1\subsetneq\cdots\subsetneq M_n=M.
\]
If there is no finite bound on \(n\), the length is infinite. See
\cite[Tag~00IU]{stacksproject}.  All module lengths used below are finite.
When \(R\) is a local \(\overline{\F}_q\)-algebra with residue field \(\overline{\F}_q\),
\(\operatorname{length}_R(M)=\dim_{\overline{\F}_q} M\) for every module \(M\)
of finite length.

For a closed singular point \(x\in C_0\), the \emph{delta invariant}
of \(C_0\) at \(x\) is
\[
 \delta_x=\operatorname{length}_{\mathcal O_{C_0,x}}
 ((\nu_*\mathcal O_C/\mathcal O_{C_0})_x).
\]
For a geometric singular point \(P\in\Sigma(C_0)\), define
\(\delta_P\) in the same way using \(\overline C_0\) and its
normalization \(\overline\nu\).
The genus formula for the normalization gives
\[
   3-g(C)
   =\delta(C_0)
   =\sum_{x\in\operatorname{Sing}(C_0)}
      [\kappa(x):\F_q] \,\delta_x
   =\sum_{P\in\Sigma(C_0)}\delta_P.
\]
See \cite[Ch.~7, Proposition~5.4 ]
{liu2002algebraic}.  Thus the normalization of a singular plane quartic has
genus two, one, or zero when its total delta invariant is one, two, or three,
respectively.

\begin{figure}[!htbp]
\centering
\begingroup
\begin{tikzpicture}[
  x=1cm,y=1cm,font=\small,
  branch/.style={draw=black,line width=0.95pt},
  rationalpoint/.style={circle,fill=black,inner sep=1.45pt},
  conjugatepoint/.style={circle,draw=black,fill=white,line width=0.65pt,inner sep=1.65pt},
  every node/.style={outer sep=2pt}]
\path[use as bounding box] (-2.45,-2.0) rectangle (12.85,2.85);

\foreach \panel/\shift in {1/0,2/5.2,3/10.4}{
\begin{scope}[xshift=\shift cm]
  \ifnum\panel=1
    \node at (0,2.6) {split node};
  \else\ifnum\panel=2
    \node at (0,2.6) {nonsplit node};
  \else
    \node at (0,2.6) {cusp};
  \fi\fi
  \node[anchor=east] at (-1.5,1.5) {\(\overline C\)};
  \draw[branch] plot[domain=-1.3:1.3,samples=80,variable=\t]
    ({\t},{1.5+0.21*sin(160*\t)});
  \pgfmathsetmacro{\normalizationLeftHeight}{1.5-0.21*sin(128)}
  \pgfmathsetmacro{\normalizationRightHeight}{1.5+0.21*sin(128)}
  \ifnum\panel=3
    \node[rationalpoint] at (0,1.5) {};
    \node at (0,0.96) {one rational preimage};
  \else
    \ifnum\panel=1
      \node[rationalpoint] at (-0.8,\normalizationLeftHeight) {};
      \node[rationalpoint] at (0.8,\normalizationRightHeight) {};
      \node at (0,0.96) {two rational preimages};
    \else
      \node[conjugatepoint] at (-0.8,\normalizationLeftHeight) {};
      \node[conjugatepoint] at (0.8,\normalizationRightHeight) {};
      \draw[<->,line width=0.5pt] (-0.8,1.53)
        .. controls (-0.6,2.35) and (0.6,2.35) .. (0.8,1.87);
      \node[fill=white,inner sep=1pt] at (0,2.21) {\(\Frob_q\)};
      \node at (0,0.96) {no rational preimage};
    \fi
  \fi
  \draw[->,line width=0.65pt] (0,0.63) -- (0,-0.12)
    node[midway,right=3pt] {\(\overline\nu\)};
  \node[anchor=east] at (-1.5,-0.72) {\(\overline C_0\)};
  \ifnum\panel=3
    \draw[branch] plot[domain=-90:90,samples=121,variable=\t]
      ({-0.6+1.65*cos(\t)^2},{-0.72+1.4*sin(\t)*cos(\t)^3}) -- cycle;
    \node[rationalpoint,label={left:\(P\)}] at (-0.6,-0.72) {};
  \else
    \draw[branch] plot[domain=0:360,samples=161,variable=\t]
      ({1.18*sin(\t)},{-0.72+0.48*sin(2*\t)}) -- cycle;
    \node[rationalpoint,label={below:\(P\)}] at (0,-0.72) {};
  \fi
  \ifnum\panel=1
    \node at (0,-1.7) {contribution: \(2-1=1\)};
  \else\ifnum\panel=2
    \node at (0,-1.7) {contribution: \(0-1=-1\)};
  \else
    \node at (0,-1.7) {contribution: \(1-1=0\)};
  \fi\fi
\end{scope}
}
\end{tikzpicture}
\endgroup
\caption{Normalization above a rational singular point \(P\).
The upper curves show smooth portions of \(\overline C\), with the
distinct preimages of \(P\) marked.  The lower curves are schematic
affine views illustrating the singularities at \(P\) over
\(\overline{\F}_q\).  Filled points are \(\F_q\)-rational,
and the two open points are exchanged by \(\Frob_q\).
The displayed differences give the contributions to \(\Delta(C_0)\).
The cusp has one geometric preimage, although its fiber has length two.}
\label{fig:normalization-trace-shift}
\end{figure}
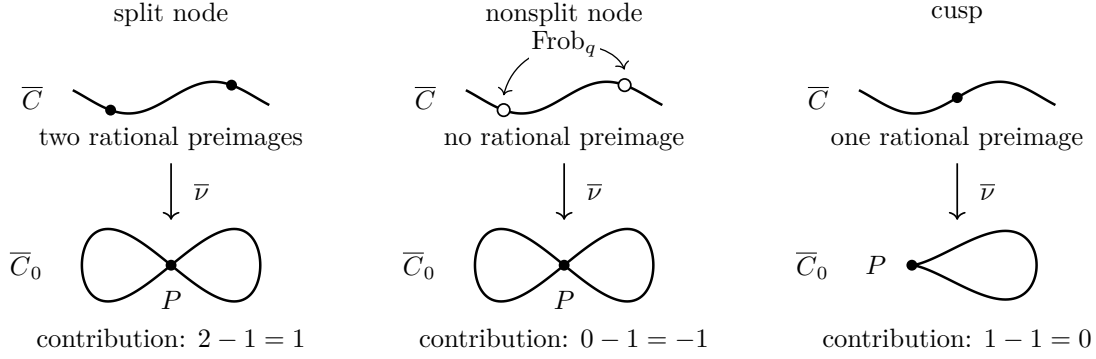

The map \(\nu\) is an isomorphism over the smooth locus.  We write
\(\mathcal B(C_0)^{\Frob_q}\) and \(\Sigma(C_0)^{\Frob_q}\) for the
subsets of \(\mathcal B(C_0)\) and \(\Sigma(C_0)\) fixed by Frobenius.
They consist of the rational branches and rational singular points,
respectively.  To count rational points on \(C_0\), we therefore remove
the rational branches from the count on \(C\) and add the rational
singular points.  This gives
\begin{equation}\label{eq:trace-shift}
   \#C_0(\F_q)
   =\#C(\F_q)-\#\mathcal B(C_0)^{\Frob_q}+\#\Sigma(C_0)^{\Frob_q}.
\end{equation}
Define
\(\Delta(C_0)=\#\mathcal B(C_0)^{\Frob_q}
-\#\Sigma(C_0)^{\Frob_q}\).  Then
\(\#C_0(\F_q)=\#C(\F_q)-\Delta(C_0)\), so the Frobenius trace of
\(C_0\) is the trace of its normalization plus \(\Delta(C_0)\).
Figure~\ref{fig:normalization-trace-shift} illustrates the contribution
of a rational node or cusp to this difference.

At each singular point \(P\), we record its multiplicity \(m(P)\)
and its delta invariant \(\delta_P\), defined above.  The multiplicity
\(m(P)\) is the smallest total degree of a nonzero term in the
defining series \(f(u,v)\) at \(P\).  Each point has its own pair of integers.
Frobenius preserves these numbers:
\(m(\Frob_q(P))=m(P)\) and
\(\delta_{\Frob_q(P)}=\delta_P\).

\Needspace{4\baselineskip}
\begin{defn}\label{def:galois-profile}
The \emph{Galois refined singularity profile} of \(C_0\) is
\[
 \tau(C_0)=
 \bigl[\mathcal B(C_0)\xrightarrow{\pi_{C_0}}\Sigma(C_0),
       (m(P),\delta_P)_{P\in\Sigma(C_0)}\bigr],
\]
where the brackets denote the isomorphism class, with the Frobenius action
on both finite sets and the pair \((m(P),\delta_P)\) attached to each
singular point \(P\).  For another quartic \(C'_0\), write
\(m'(P')\) and \(\delta'_{P'}\) for the two integers at a singular
point \(P'\).  The two quartics have the same profile
exactly when there are Frobenius equivariant bijections
\(\alpha\colon\mathcal B(C_0)\xrightarrow{\sim}\mathcal B(C'_0)\)
and \(\beta\colon\Sigma(C_0)\xrightarrow{\sim}\Sigma(C'_0)\)
that preserve the map to the singular points,
\[
 \pi_{C'_0}\circ\alpha=\beta\circ\pi_{C_0},
\]
and the two integers at corresponding points, so that for every
\(P\in\Sigma(C_0)\),
\[
\begin{aligned}
 m'(\beta(P))&=m(P),\\
 \delta'_{\beta(P)}&=\delta_P.
\end{aligned}
\]
\end{defn}
\par

The map \(\pi_{C_0}\) already determines the branch number \(r(P)\).
For the quartic singularity types considered here, the triples
\((m(P),\delta_P,r(P))\) in Table~\ref{tab:local-delta-classification}
are all distinct. See \cite[\S9]{cox2013singularities}.
Thus the profile determines the local type at each point.  Conversely,
each type determines these three integers.  We use the labels in the
table as names for the types recovered from the numerical data.
The sum of the delta invariants at the singular points determines the
genus of the normalization.

For example, a node and a cusp both have \((m(P),\delta_P)=(2,1)\),
but the fibers over them have two and one branches, respectively.
If these are the only two singular points, Frobenius must fix each point
because it preserves the branch number.  Two nodes may instead form a
conjugate pair.  At a rational node, Frobenius fixes both branches in the
split case and exchanges them in the nonsplit case.
These possibilities are all recorded by the finite sets and their map.

For a profile \(\tau\), write
\(\mathcal B_\tau\to\Sigma_\tau\) for its underlying map of finite
Frobenius sets, and define
\(\Delta_\tau=\#\mathcal B_\tau^{\Frob_q}
-\#\Sigma_\tau^{\Frob_q}\).
This number depends only on the isomorphism class of \(\tau\), and every
quartic of profile \(\tau\) satisfies
\begin{equation}\label{eq:profile-point-shift}
   \#C_0(\F_q)=\#C(\F_q)-\Delta_\tau.
\end{equation}

We recall the multiplicities of branches, their multiplicity sequences,
and the intersection multiplicities between distinct branches.
We use these to distinguish the singularity types in
Table~\ref{tab:local-delta-classification}.

We describe these local invariants over \(\overline{\F}_q\).
At a singular point \(P\), choose
coordinates \(u,v\) centered at \(P\), as above.  For a branch \(B\)
at \(P\), let \(Q\in\overline C\) be its corresponding point on the
normalization.  Choose a local parameter \(t\) at \(Q\), meaning a
function that vanishes to order one there.  The completed local ring at
\(Q\) is then \(\overline{\F}_q[[t]]\) \cite[Theorem~16.9]{kunz2005plane}.
Pulling back \(u,v\) by \(\overline\nu\) gives formal power series
\(u(t),v(t)\in t\overline{\F}_q[[t]]\) satisfying \(f_B(u(t),v(t))=0\).
This gives a formal parametrization of~\(B\).
For a nonzero power series,
\(\operatorname{ord}_t\) is the exponent of its first nonzero term.
We set \(\operatorname{ord}_t(0)=\infty\).  The multiplicity of \(B\) is
\[
 m_0(B)=\min\{\operatorname{ord}_t u(t),
                   \operatorname{ord}_t v(t)\}.
\]
In particular, the branch is smooth exactly when its multiplicity is one
\cite[Remark~2.4(a)]{nguyen2016invariants}.
The point multiplicity is
\(m(P)=\sum_{B\in\pi_{C_0}^{-1}(P)}m_0(B)\), the sum over its
branches.  It equals \(r(P)\) exactly when all branches at \(P\) are smooth.

For a finite scheme \(Z\) over \(\overline{\F}_q\), its
\emph{length} is
\[
 \operatorname{length}(Z)=\sum_{z\in Z(\overline{\F}_q)}\dim_{\overline{\F}_q}\mathcal O_{Z,z},
\]
where \(\mathcal O_{Z,z}\) is the local ring at \(z\) and
\(\dim_{\overline{\F}_q}\) denotes vector space dimension.  We use this for finite
intersections and fibers.  For example, the doubled point \(2Q\) on
the smooth curve \(\overline C\) has local ring \(\overline{\F}_q[[t]]/(t^2)\).
Its basis over \(\overline{\F}_q\) is \(1,t\), so it has length two although it
consists of only one point.
\par

These multiplicities also describe the fiber of the normalization over
\(P\).  The local ring of this fiber at \(Q\) is
\(\overline{\F}_q[[t]]/(u(t),v(t))\).  A nonzero power series is a power of \(t\)
times a unit, so the ideal \((u(t),v(t))\) is \((t^{m_0(B)})\).
Thus the length of the fiber at \(Q\) is
\begin{equation}\label{eq:normalization-fiber-multiplicity}
 \dim_{\overline{\F}_q} \overline{\F}_q[[t]]/(u(t),v(t))
 =\min\{\operatorname{ord}_t u(t),\operatorname{ord}_t v(t)\}
 =m_0(B).
\end{equation}
Summing over the branches gives
\(\operatorname{length}(\overline\nu^{-1}(P))=m(P)\).
Here length counts each point with its multiplicity, whereas \(r(P)\)
counts the distinct points in the fiber.  At a smooth point, the
normalization is an isomorphism and its fiber has length one.

The \emph{multiplicity sequence} records what happens under successive
blowups.  After blowing up the singular point, follow the strict transform
of \(B\), namely the closure of the part of the branch away from that point,
and record its multiplicity at the point above it.  Repeat until the branch
first becomes smooth.  Let \(N_B\) be the number of blowups performed,
and let \(m_i(B)\) be the multiplicity after \(i\) blowups, for
\(0\leq i\leq N_B\).  We write the resulting list as
\[
 \mathbf m(B)=(m_0(B):m_1(B):\cdots:m_{N_B}(B)).
\]
The last entry is \(m_{N_B}(B)=1\).  If \(B\) is already smooth, then \(N_B=0\)
and its sequence is \((1)\).

This sequence can be computed directly from a parametrization.  If
\(\operatorname{ord}_t u\leq\operatorname{ord}_t v\), the relevant blowup
chart has coordinates \(u\) and \(v/u-c\), where \(c\) is
the constant term of \(v/u\).  For example, the branch
\((u,v)=(t^2,t^5)\) becomes \((t^2,t^3)\) after one blowup and
\((t^2,t)\) after the next.  Thus \(N_B=2\), its multiplicity sequence
is \((2:2:1)\), and its delta invariant is two.

More generally, the blowup formula gives
\begin{equation}\label{eq:branch-delta-multiplicities}
 \delta(B)=\sum_{i=0}^{N_B}\binom{m_i(B)}2.
\end{equation}
See \cite[Remark~2.10]{nguyen2016invariants}.
\par

For two distinct branches \(B,B'\) at the same point, choose their
irreducible defining series \(f_B,f_{B'}\in \overline{\F}_q[[u,v]]\), as above.
Their \emph{intersection multiplicity} is
\[
 I(B,B')=\dim_{\overline{\F}_q} \overline{\F}_q[[u,v]]/(f_B,f_{B'}).
\]
\par
Two smooth branches with different tangent lines have intersection
multiplicity one.  The branches \(v=0\) and \(v=u^2\) have intersection
multiplicity two.  Thus \(\mathbf m(B)\) describes each branch on its own,
whereas \(I(B,B')\) measures the contact between different branches.
If \(B_1,\ldots,B_s\) are the branches at \(P\), then
\begin{equation}\label{eq:local-delta-branches}
 \delta_P=\sum_{i=1}^s\delta(B_i)
             +\sum_{1\leq i<j\leq s}I(B_i,B_j)
\end{equation}
\cite[Remark~2.4(d)]{nguyen2016invariants}.  In particular, singularities
with two smooth branches can have the same branch multiplicity sequences
but different delta invariants because the intersection multiplicities
differ.

Since a plane quartic has arithmetic genus three, every singular point
satisfies \(\delta_P\leq3\).  The standard formulas for the delta invariant
give the branch configurations in Table~\ref{tab:local-delta-classification}.
We use the labels of \cite[\S9]{cox2013singularities}, expressed in terms
of the branch data above.  As throughout the paper, we assume
\(\operatorname{char}(\F_q)>3\).

\begin{table}[!htbp]
\centering
\caption{Local singularity configurations for absolutely irreducible plane quartics.}
\label{tab:local-delta-classification}
\renewcommand{\arraystretch}{1.14}
\begin{tabular}{@{}cc>{\raggedright\arraybackslash}p{0.64\textwidth}@{}}
\toprule
Label & \((m(P),\delta_P,r(P))\) & Intrinsic branch data \\
\midrule
\(A_1\) & \((2,1,2)\) & two smooth branches with intersection multiplicity one \\
\(A_2\) & \((2,1,1)\) & one branch with multiplicity sequence \((2:1)\) \\
\(A_3\) & \((2,2,2)\) & two smooth branches with intersection multiplicity two \\
\(A_4\) & \((2,2,1)\) & one branch with multiplicity sequence \((2:2:1)\) \\
\(D_4\) & \((3,3,3)\) & three smooth branches, pairwise of intersection multiplicity one \\
\(D_5\) & \((3,3,2)\) & one smooth branch and one \((2:1)\) branch, of intersection multiplicity two \\
\(A_5\) & \((2,3,2)\) & two smooth branches with intersection multiplicity three \\
\(E_6\) & \((3,3,1)\) & one branch with multiplicity sequence \((3:1)\) \\
\(A_6\) & \((2,3,1)\) & one branch with multiplicity sequence \((2:2:2:1)\) \\
\bottomrule
\end{tabular}
\end{table}
For comparison with the literature, we also recall the \emph{value
semigroup} of a branch.  Its completed local ring is the subring
\(R_B=\overline{\F}_q[[u(t),v(t)]]\subset \overline{\F}_q[[t]]\), and its value semigroup is
\[
 \Gamma(B)=\{\operatorname{ord}_t h:0\ne h\in R_B\}
             \subset\mathbb Z_{\geq0}.
\]
It contains zero and is closed under addition because orders add under
multiplication.  Its complement is finite, and the number of missing
nonnegative integers is \(\delta(B)\).
See~\mbox{\cite[\S\S4.3--4.4]{campillo1980algebroid}}.  Here
\mbox{\(\langle a,b\rangle=\{ia+jb:i,j\in\mathbb Z_{\geq0}\}\)}.
For example, \(\overline{\F}_q[[t^2,t^7]]\) has value semigroup
\(\langle2,7\rangle=\{0,2,4,6,7,8,\ldots\}\).  Its missing values
are \(1,3,5\), so its delta invariant is three.  The four unibranch types
\(A_2,A_4,E_6,A_6\) in the table have value semigroups
\(\langle2,3\rangle\), \(\langle2,5\rangle\),
\(\langle3,4\rangle\), and \(\langle2,7\rangle\), respectively.
These give another description of the listed branch types.  We will check
the types needed in the direct counts by vanishing orders and blowups,
without a separate classification of value semigroups.

\begin{remark}
\label{rem:positive-characteristic-labels}
The symbols \(A_i,D_i,E_6\) refer only to the branch data in
Table~\ref{tab:local-delta-classification}.  We do not use them to claim
that a defining series in \(\overline{\F}_q[[u,v]]\) is right equivalent to a characteristic zero
ADE normal form, or that the singularity is right simple.  These properties
can change in positive characteristic.  For example, in characteristic
\(p>2\), the right simple ranges are
\(A_k\) with \(k\le p-2\) and \(D_k\) with \(k<p\)
\cite[Theorem~3.2]{greuelnguyen2016rightsimple}.  Thus, for instance,
\(A_4,A_5,D_5\) in characteristic five and \(A_6\) in characteristic seven
do not have the characteristic zero right simple meaning.

The local calculations below use vanishing orders, multiplicity sequences,
and branch intersection multiplicities.  The point counts also use the
Frobenius action on branches.  Bergstr\"om's
counts used in genus zero come from his odd characteristic formulas,
and the direct local computations require only that two and three be
invertible.  Thus the assumption \(\operatorname{char}(\F_q)>3\) is enough.
For all the types in Table~\ref{tab:local-delta-classification} to have
the characteristic zero right simple meaning, it is enough to assume
\(\operatorname{char}(\F_q)>7\).
\end{remark}

We next define what we count on a fixed normalization.  Let \(C/\F_q\)
be a smooth projective geometrically connected curve.  A base point free linear series
\(g^2_4\) on \(C\) is an isomorphism class of pairs \((L,V)\), where \(L\)
is a line bundle of degree four on \(C\) and \(V\subset H^0(C,L)\) is a
three dimensional base point free subspace.  Two pairs \((L,V)\) and
\((L',V')\) are
isomorphic if an isomorphism \(L\simeq L'\) carries \(V\) onto \(V'\).
The pair defines a morphism
\(\phi_{L,V}\colon C\longrightarrow\PP(V^\vee)\simeq\PP^2\).
Let \(\mathcal G^2_4(C)\) be the set of such linear series for which
\(\phi_{L,V}\) is geometrically birational onto its image.  For
\(\gamma=[(L,V)]\in\mathcal G^2_4(C)\), put
\(C_\gamma=\phi_{L,V}(C)\subset\PP(V^\vee)\).
For a profile \(\tau\), let \(\mathcal G^2_{4,\tau}(C)\) be the subset of
\(\mathcal G^2_4(C)\) consisting of the \(\gamma\) for which
\(\tau(C_\gamma)=\tau\), and write
\(m_\tau(C)=\#\mathcal G^2_{4,\tau}(C)\).  Thus \(m_\tau(C)\) counts
plane linear series on the fixed curve \(C\).  We do not identify series
that differ by an element of \(\Aut_{\F_q}(C)\).

\begin{prop}\label{prop:irr-strata-contrib}
Let \(\tau\) be a profile whose quartics have normalizations of genus \(g\).
Its contribution to the weight enumerator is
\[
  W_{\cC_{2,4}}^{\Irr,\tau}(X,Y)
  =(q-1)\,\#\PGL_3(\F_q)
  \sum_{[C]}
  \frac{m_\tau(C)
  X^{\#C(\F_q)-\Delta_\tau}
  Y^{q^2+q+1-\#C(\F_q)+\Delta_\tau}}
  {\#\Aut_{\F_q}(C)}.
\]
The sum is over representatives of the \(\F_q\)-isomorphism classes of
smooth projective genus \(g\) curves.  Here \(\Aut_{\F_q}(C)\)
is the full automorphism group of the underlying curve.
\end{prop}

\begin{proof}
Fix \(C\).  For each
\(\gamma=[(L,V)]\in\mathcal G^2_{4,\tau}(C)\), the set of projective
isomorphisms \(\PP(V^\vee)\simeq\PP^2_{\F_q}\)
has \(\#\PGL_3(\F_q)\) elements.  Hence the \(m_\tau(C)\) linear series
give \(\#\PGL_3(\F_q)m_\tau(C)\) normalization maps \(C\to\PP^2\).  On
the other hand, an embedded quartic
\(C_0\subset\PP^2\) whose normalization is \(\F_q\)-isomorphic to \(C\)
has exactly \(\#\Aut_{\F_q}(C)\) normalization maps from \(C\).  After
choosing one map, we obtain all the others uniquely by composing it with an
automorphism of \(C\).  Therefore the number of embedded quartics of
profile \(\tau\) with normalization isomorphic to \(C\) is
\(\#\PGL_3(\F_q)m_\tau(C)/\#\Aut_{\F_q}(C)\).
Each embedded quartic has \(q-1\) nonzero scalar defining equations, and
\eqref{eq:profile-point-shift} gives its number of rational points.  Summing
over \([C]\) proves the stated formula.
\end{proof}

\subsection{Geometric genus two}

Fix a smooth projective genus two curve \(C/\F_q\).  For
\(M\in(\Sym^2 C)(\F_q)\), put \(\mathcal L_M=\omega_C(M)\), and let
\(\phi_M=\phi_{|\mathcal L_M|}\colon C\longrightarrow\PP^2\).
By Riemann--Roch and \cite[Corollary~2.19]{eisenbud2024practice},
every line bundle of degree four on \(C\) has three global sections and is
base point free.  The classification of genus two plane linear
series in \cite[\S6.3]{eisenbud2024practice} gives the following description.
If \(M\in|K_C|\), then
\(\mathcal L_M\simeq\omega_C^{\otimes2}\), and \(\phi_M\) is the canonical
hyperelliptic double cover followed by the quadratic Veronese embedding of
\(\PP^1\).  If \(M\notin|K_C|\), define
\(C_M:=\phi_M(C)\subset\PP^2\)
and let \(\nu_M\colon C\to C_M\) be the induced map.  Then \(C_M\) is an
absolutely irreducible plane quartic, \(\nu_M\) is its normalization, and the
scheme theoretic fiber of \(\nu_M\) over the unique singular point of \(C_M\)
is \(M\).  Thus, after base change to \(\overline{\F}_q\) and writing
\(M=P+Q\), the singularity is a node when \(P\ne Q\) and a cusp when
\(P=Q\).

Conversely, let
\(
\varphi \colon C \longrightarrow \mathbf P^2
\)
be birational onto a plane quartic, and put
\(
\mathcal L=\varphi^*\mathcal O_{\mathbf P^2}(1).
\)
Then $\deg \mathcal L=4$, and Riemann--Roch gives
$h^0(C,\mathcal L)=3$. Since the image of $\varphi$ is not contained
in a line, the pullbacks of the linear forms on $\mathbf P^2$ span
all of $H^0(C,\mathcal L)$. Hence $\varphi$ is defined by the complete
linear series $|\mathcal L|$. Birationality excludes
$\mathcal L\simeq \omega_C^{\otimes 2}$, since the latter gives the
hyperelliptic double cover followed by the quadratic Veronese
embedding.

Now $\mathcal L\otimes\omega_C^{-1}$ has degree two, and
Riemann--Roch gives
\(
h^0\!\left(C,\mathcal L\otimes\omega_C^{-1}\right)\geq 1.
\)
Thus
\(
\mathcal L\otimes\omega_C^{-1}\simeq \mathcal O_C(M)
\)
for some $M\in(\operatorname{Sym}^2 C)(\mathbf F_q)$. Since
$\mathcal L\not\simeq\omega_C^{\otimes 2}$, we have
$M\notin |K_C|$. Moreover, a degree two line bundle on a genus two
curve has more than one section only when it is isomorphic to
$\omega_C$. Hence $M$ is unique. Therefore
$\mathcal L\simeq\omega_C(M)$, and $\varphi$ agrees with $\phi_M$
up to projective coordinates. Consequently,
\(
(\operatorname{Sym}^2 C)(\mathbf F_q)
   \setminus |K_C|(\mathbf F_q)
\)
parametrizes the elements of $\mathcal G^2_4(C)$,
equivalently the birational plane quartic maps
$C\to\mathbf P^2$ up to projective changes of coordinates in the target.

Since \(C_M\) has a unique geometric singular point, that point is
\(\F_q\)-rational.  The following table lists the three possible Frobenius
types of \(M\) and their trace shifts \(\Delta_M=\Delta(C_M)\).
In every row we require \(M\notin|K_C|\).
\begin{center}
\renewcommand{\arraystretch}{1.2}
\begin{tabular}{@{}T{5.1cm}T{3.6cm}cc@{}}
\toprule
\(M\) & \(\F_q\)-profile & trace shift \(\Delta_M\) & \(\#C_M(\F_q)\) \\
\midrule
\(P+Q\), with distinct \(P,Q\in C(\F_q)\)
& split node & 1 & \(\#C(\F_q)-1\)\\
\(2P\), with \(P\in C(\F_q)\)
& rational cusp & 0 & \(\#C(\F_q)\)\\
\(P+\Frob_q(P)\), with \(P\in C(\F_{q^2})\setminus C(\F_q)\)
& nonsplit node & \(-1\) & \(\#C(\F_q)+1\)\\
\bottomrule
\end{tabular}
\end{center}
The last column follows from
\eqref{eq:trace-shift}.  Denote the corresponding multiplicities on the
fixed normalization by \(m_{\mathrm{split}}(C)\),
\(m_{\mathrm{cusp}}(C)\), and \(m_{\mathrm{nonsplit}}(C)\).

Let $h\colon C\to\mathbf P^1$ be the canonical hyperelliptic double cover,
and write $\operatorname{Weier}(C)\subset C$ for its ramification locus,
equivalently the set of Weierstrass points of $C$. Write $J(C)=\operatorname{Pic}^0(C)$ for the Jacobian of $C$.

\begin{prop}\label{prop:g2-fixed-counts}
The split node multiplicity is
\[
 m_{\mathrm{split}}(C)
 =\frac{\#C(\F_q)(\#C(\F_q)-2)
 +\#\operatorname{Weier}(C)(\F_q)}2.
\]
The cusp multiplicity is
\[
 m_{\mathrm{cusp}}(C)
 =\#C(\F_q)-\#\operatorname{Weier}(C)(\F_q),
\]
and the nonsplit node multiplicity is
\[
 m_{\mathrm{nonsplit}}(C)
 =\frac{\#C(\F_{q^2})+\#\operatorname{Weier}(C)(\F_q)-2q-2}2.
\]
Their sum is
\[
   m_{\mathrm{split}}(C)+m_{\mathrm{cusp}}(C)
   +m_{\mathrm{nonsplit}}(C)
   =
   \frac{\#C(\F_q)^2+\#C(\F_{q^2})}{2}-(q+1)
   =
   \#J(C)(\F_q)-1.
\]
\end{prop}

\begin{proof}
The
\(q+1\) rational fibers of \(h\) include
\(\#\operatorname{Weier}(C)(\F_q)\) ramified fibers.  Among the remaining
fibers,
\((\#C(\F_q)-\#\operatorname{Weier}(C)(\F_q))/2\) are split and
nonramified.  The number of nonsplit fibers is therefore
\[
 q+1-\#\operatorname{Weier}(C)(\F_q)
 -\frac{\#C(\F_q)-\#\operatorname{Weier}(C)(\F_q)}2.
\]
On the other hand, there are \(\binom{\#C(\F_q)}2\) rational effective
degree two divisors of split type.  There are \(\#C(\F_q)\) double divisors,
and there are
\((\#C(\F_{q^2})-\#C(\F_q))/2\) nonsplit divisors.
The excluded divisors are exactly the fibers of \(h\).  Subtracting the
split nonramified fibers from the split divisors gives
\[
 \binom{\#C(\F_q)}2
 -\frac{\#C(\F_q)-\#\operatorname{Weier}(C)(\F_q)}2
 =m_{\mathrm{split}}(C).
\]
Removing the ramified fibers from the double divisors leaves
\(\#C(\F_q)-\#\operatorname{Weier}(C)(\F_q)\), which is
\(m_{\mathrm{cusp}}(C)\).  Finally, subtracting the nonsplit canonical
fibers from the nonsplit divisors gives
\[
 \frac{\#C(\F_{q^2})-\#C(\F_q)}2
 -\left(q+1-\#\operatorname{Weier}(C)(\F_q)
 -\frac{\#C(\F_q)-\#\operatorname{Weier}(C)(\F_q)}2\right),
\]
which simplifies to \(m_{\mathrm{nonsplit}}(C)\).  Adding these identities
gives the stated total.  To prove the final equality, recall that the
classification above identifies all birational degree four plane models of
\(C\), up to projective coordinates, with
\(\Pic^4(C)(\F_q)\setminus\{\omega_C^{\otimes2}\}\).
Translation by \((\omega_C^{\otimes2})^{-1}\) identifies \(\Pic^4(C)\) with
\(J(C)\), so this set has \(\#J(C)(\F_q)-1\) elements.
\end{proof}

Combining Propositions~\ref{prop:g2-fixed-counts} and
\ref{prop:irr-strata-contrib} gives
\begin{align}
  W_{\cC_{2,4}}^{\Irr,\delta=1}(X,Y)
  &=(q-1)\,\#\PGL_3(\F_q)
  \sum_{[C]\in\mathcal M_2(\F_q)}
  \frac{1}{\#\Aut_{\F_q}(C)}
  \bigg(
      m_{\mathrm{split}}(C)
      X^{\#C(\F_q)-1}Y^{q^2+q+2-\#C(\F_q)}
  \nonumber\\*
  &{}+m_{\mathrm{cusp}}(C)
      X^{\#C(\F_q)}Y^{q^2+q+1-\#C(\F_q)}
      +m_{\mathrm{nonsplit}}(C)
      X^{\#C(\F_q)+1}Y^{q^2+q-\#C(\F_q)}
  \bigg).
\label{eq:g2-weight-contribution}
\end{align}

\subsection{Geometric genus one}
\label{subsec:geometric-genus-one}

Let \(E/\F_q\) be a smooth projective curve of genus one and \(L\in\Pic^4(E)(\F_q)\).  The complete linear series \(|L|\) embeds
\(E\) as the elliptic normal quartic
\(\iota_L\colon E\hookrightarrow
\PP(H^0(E,L)^\vee)\simeq\PP^3\).
Its image is the complete intersection of two independent members of a pencil
\(\Lambda_L\simeq\PP^1\) of quadrics
\cite[\S4.5]{eisenbud2024practice}.  Here a pencil means the family
\(V(aF_0+bF_1)\), with \([a:b]\in\PP^1\), generated by two independent
quadratic equations \(F_0,F_1\) vanishing on \(E\).  A three dimensional
subspace \(V\subset H^0(E,L)\) determines a unique projection center
\(x\in\PP(H^0(E,L)^\vee)\), namely the point represented by the
linear functionals that vanish on \(V\).  The series is base point free
exactly when \(x\notin\iota_L(E)\). Thus, on the fixed curve \(E\), we count pairs \((L,x)\) that give
birational projections.  The quadrics describe which singularities a center
produces.  We then use effective divisors on \(E\) to count these centers. A ruling line, or generator, means a line contained in a quadric.  Over
\(\overline{\F}_q\), a smooth quadric has two families of such lines,
with one line from each family through every point.  A quadric cone is
the union of the lines joining its vertex to a smooth conic.  Away from
the vertex, each point lies on just one generator.

For an effective divisor \(M\), its support \(\operatorname{Supp}(M)\)
is the set of geometric points that occur in \(M\), without
multiplicities.  Two divisors have disjoint support if these sets have no
common point.  For \(M\) of degree two, \(\langle M\rangle\) denotes its
secant line when \(M\) is reduced and its tangent line when \(M=2P\).
We write \([M]\) for its linear equivalence class in
\(\Pic^2(E_{\overline{\F}_q})\).
The complete linear series \(|M|\) consists of the effective divisors
linearly equivalent to \(M\).  Riemann--Roch shows that it is a base point
free pencil isomorphic to \(\PP^1\).

\begin{prop}
\label{prop:g1-projection-geometry}
Let \(x\in\PP^3(\F_q)\setminus\iota_L(E)(\F_q)\), let
\(\pi_x\colon\PP^3\dashrightarrow\PP^2\) be projection from \(x\), and
write \(\phi_x=\pi_x\circ\iota_L\) and \(C_x=\phi_x(E)\).
There is a unique quadric \(Q_x\in\Lambda_L\) containing \(x\).
Over \(\overline{\F}_q\), the three cases have the following divisor
descriptions.
\begin{enumerate}[label=\textup{(\roman*)}]
\item If \(Q_x\) is smooth, its two ruling lines through \(x\) meet
      \(E\) in disjoint effective divisors \(M_1,M_2\) of degree two,
      with \(L\simeq\mathcal O_E(M_1+M_2)\) and \([M_1]\neq[M_2]\).
      The map \(\phi_x\) is birational onto a quartic with two
      singularities, each with delta invariant one.  The singularity
      corresponding to \(M_i\) is a node when \(M_i\) is reduced and a
      cusp when \(M_i=2P\).

\item If \(Q_x\) is a cone with vertex \(v\neq x\), the line joining
      \(v\) and \(x\) meets \(E\) in an effective divisor \(M\) of
      degree two, with \(L\simeq\mathcal O_E(2M)\).
      The map \(\phi_x\) is birational onto a quartic with one
      singularity of delta invariant two.  It is an \(A_3\) singularity
      when \(M\) is reduced and an \(A_4\) singularity when \(M=2P\).

\item If \(x\) is the vertex of \(Q_x\), then \(\phi_x\) has degree
      two onto a conic.  Its fibers form a degree two pencil \(|M|\),
      with \(L\simeq\mathcal O_E(2M)\).
\end{enumerate}
Conversely, two effective degree two divisors \(M_1,M_2\) with disjoint
support and \(L\simeq\mathcal O_E(M_1+M_2)\) determine the center
\(x=\langle M_1\rangle\cap\langle M_2\rangle\).
This gives case (i) when \([M_1]\neq[M_2]\) and case (iii) when
\([M_1]=[M_2]\).
For an effective degree two divisor \(M\) with
\(L\simeq\mathcal O_E(2M)\), the lines \(\langle M'\rangle\), for
\(M'\in|M|\), form a quadric cone \(Q\) with vertex \(v\).
The centers in \(\langle M\rangle\) outside
\(\operatorname{Supp}(M)\cup\{v\}\) give precisely the projections
in case (ii) with singular fiber \(M\).
\end{prop}

\begin{figure}[!htbp]
\centering
\begingroup
\pgfmathdeclarefunction{quarticSaddleR}{1}{%
  \pgfmathparse{sqrt(2/(1+sqrt(1+sin(2*#1)^2)))}}
\pgfmathdeclarefunction{quarticConeR}{2}{%
  \pgfmathparse{(1.2-1.2*sin(#1)+#2*sqrt(max(0,(1.2-1.2*sin(#1))^2-4.14)))/2}}
\pgfmathsetmacro{\quarticAngle}{acos((sqrt(4.14)-1.2)/1.2)}
\pgfmathsetmacro{\quarticStart}{270-\quarticAngle}
\pgfmathsetmacro{\quarticEnd}{270+\quarticAngle}
\begin{tikzpicture}[
  x=1cm,y=1cm,font=\small,
  surface/.style={draw=black!35,line width=0.3pt},
  curve/.style={draw=black,line width=1.05pt},
  ruling/.style={draw=black,line width=0.65pt},
  fiberpoint/.style={circle,fill=black,inner sep=1.35pt},
  projectioncenter/.style={circle,draw=black,fill=white,line width=0.6pt,inner sep=1.65pt},
  every node/.style={outer sep=2pt}]
\path[use as bounding box] (-2.45,-3.95) rectangle (12.85,4.7);

\begin{scope}
  \node at (0,4.4) {\textup{(i)} Smooth \(Q_x\)};
  \fill[black!3]
    (-2.0125,2.405125) -- (0.6325,0.754875) --
    (2.0125,3.785125) -- (-0.6325,2.674875) -- cycle;
  \foreach \a in {-1.15,-0.77,-0.38,0.38,0.77,1.15}{
    \draw[surface] plot[domain=-1.15:1.15,samples=2,variable=\t]
      ({1.15*\a+0.6*\t},{2.5-0.2*\a+0.8*\t+0.45*\a*\t});
    \draw[surface] plot[domain=-1.15:1.15,samples=2,variable=\t]
      ({1.15*\t+0.6*\a},{2.5-0.2*\t+0.8*\a+0.45*\t*\a});
  }
  \draw[ruling] (-1.5525,2.77) -- (1.5525,2.23);
  \draw[ruling] (-0.81,1.42) -- (0.81,3.58);
  \draw[curve] plot[domain=0:360,samples=150,variable=\t]
    ({quarticSaddleR(\t)*(1.15*cos(\t)+0.6*sin(\t))},
     {2.5+quarticSaddleR(\t)*(-0.2*cos(\t)+0.8*sin(\t))
      +0.45*quarticSaddleR(\t)^2*cos(\t)*sin(\t)}) -- cycle;
  \foreach \p in {(1.15,2.3),(-1.15,2.7),(0.6,3.3),(-0.6,1.7)}
    \node[fiberpoint] at \p {};
  \node[projectioncenter,label={[label distance=1pt]below right:\(x\)}] at (0,2.5) {};
  \node[anchor=south east] at (-1.48,2.8) {\(\langle M_1\rangle\)};
  \node[anchor=south] at (0.79,3.61) {\(\langle M_2\rangle\)};
  \node at (-1.27,1.65) {\(E\)};
  \node at (1.85,3.98) {\(Q_x\)};
  \node at (0,0.4) {\(E\cap\langle M_i\rangle=M_i\)};
\end{scope}

\foreach \panel/\shift in {2/5.2,3/10.4}{
\begin{scope}[xshift=\shift cm]
  \ifnum\panel=2
    \node at (0,4.4) {\textup{(ii)} Cone, \(x\neq v\)};
  \else
    \node at (0,4.4) {\textup{(iii)} Cone, \(x=v\)};
  \fi
  \fill[black!3] (0,3.6) -- (-1.8675,1.35)
    arc[start angle=180,end angle=360,x radius=1.8675,y radius=0.5175] -- cycle;
  \draw[surface,dashed] (1.8675,1.35)
    arc[start angle=0,end angle=180,x radius=1.8675,y radius=0.5175];
  \draw[surface] (-1.8675,1.35) -- (0,3.6) -- (1.8675,1.35);
  \foreach \a in {215,235,250,270,290,305,325}{
    \draw[surface] (0,3.6) -- ({1.8675*cos(\a)},{1.35+0.5175*sin(\a)});
  }
  \ifnum\panel=2
    \draw[ruling] (0,3.6) -- (0,0.8325);
  \else
    \foreach \a in {250,290}{
      \draw[ruling,->] (0,3.6) -- ({1.8675*cos(\a)},{1.35+0.5175*sin(\a)});
    }
  \fi
  \draw[curve]
    plot[domain=\quarticStart:\quarticEnd,samples=100,variable=\t]
      ({0.83*quarticConeR(\t,1)*cos(\t)},
       {3.6-quarticConeR(\t,1)+0.23*quarticConeR(\t,1)*sin(\t)})
    -- plot[domain=\quarticEnd:\quarticStart,samples=100,variable=\t]
      ({0.83*quarticConeR(\t,-1)*cos(\t)},
       {3.6-quarticConeR(\t,-1)+0.23*quarticConeR(\t,-1)*sin(\t)}) -- cycle;
  \draw[surface] (-1.8675,1.35)
    arc[start angle=180,end angle=360,x radius=1.8675,y radius=0.5175];
  \node at (-1.08,2.26) {\(E\)};
  \node at (1.32,3.12) {\(Q_x\)};
  \ifnum\panel=2
    \foreach \sign in {-1,1}{
      \node[fiberpoint] at (0,{3.6-1.23*quarticConeR(270,\sign)}) {};
    }
    \node[fiberpoint,label={above left:\(v\)}] at (0,3.6) {};
    \node[projectioncenter,label={right:\(x\)}] at (0,2.124) {};
    \node[anchor=west] at (0.09,1.02) {\(\langle M\rangle\)};
    \node at (0,0.4) {\(E\cap\langle M\rangle=M\)};
  \else
    \foreach \a in {250,290}{
      \foreach \sign in {-1,1}{
        \node[fiberpoint] at ({0.83*quarticConeR(\a,\sign)*cos(\a)},
          {3.6-quarticConeR(\a,\sign)+0.23*quarticConeR(\a,\sign)*sin(\a)}) {};
      }
    }
    \node[projectioncenter,label={above:\(x=v\)}] at (0,3.6) {};
    \node at (0,0.4) {generators through \(v\)};
  \fi
\end{scope}
}

\foreach \shift in {0,5.2,10.4}{
  \draw[->,line width=0.7pt] (\shift,0.05) -- (\shift,-0.92)
    node[midway,right=3pt] {\(\phi_x\)};
  \node at (\shift,-1.24) {\(C_x\subset\PP^2\)};
}

\begin{scope}[yshift=-2.2cm]
  \foreach \sign in {-1,1}{
    \foreach \quarticLeft/\quarticRight in {-1.75/-1,-1/1,1/1.75}{
      \draw[curve] plot[domain=\quarticLeft:\quarticRight,samples=80,variable=\a]
        ({1.1*\a},{0.9*\sign*sqrt(max(0,-1-(\a)^2+sqrt(2+2*(\a)^4)))});
    }
  }
  \node[fiberpoint] at (-1.1,0) {};
  \node[fiberpoint] at (1.1,0) {};
  \node at (-1.1,-0.86) {\(\phi_x(M_2)\)};
  \node at (1.1,-0.86) {\(\phi_x(M_1)\)};
  \node at (0,-1.43) {quartic with two nodes};
\end{scope}

\begin{scope}[xshift=5.2cm,yshift=-2.2cm]
  \draw[surface,dashed] (-1.1,0) -- (1.1,0);
  \foreach \sign in {-1,1}{
    \draw[curve] plot[domain=-0.6:0.6,samples=121,variable=\a]
      ({1.15*\a*sqrt(32*(\a)^2-9*\sign+4*sqrt(82*(\a)^4-36*\sign*(\a)^2+18))/3},
       {1.75*\sign*(\a)^2});
  }
  \node[fiberpoint] at (0,0) {};
  \node at (0,-0.86) {\(\phi_x(M)\)};
  \node at (0,-1.43) {quartic with a tacnode};
\end{scope}

\begin{scope}[xshift=10.4cm,yshift=-2.2cm]
  \draw[curve] (0,0) ellipse[x radius=1.45,y radius=0.6];
  \foreach \a in {250,290}{
    \node[fiberpoint] at ({1.45*cos(\a)},{0.6*sin(\a)}) {};
  }
  \node at (0,-0.86) {fibers form \(|M|\)};
  \node at (0,-1.43) {conic, degree two map};
\end{scope}
\end{tikzpicture}
\endgroup
\caption{Projection from \(x\) in the three cases of
Proposition~\ref{prop:g1-projection-geometry}, shown schematically over
\(\overline{\F}_q\).  The marked pairs on \(E\) map to the marked points
on \(C_x\).  The reduced divisors shown give two nodes in (i) and a tacnode
in (ii).  In (iii), every fiber has length two and the image is a conic.
Only affine portions of the quartic images are drawn.  Doubled divisors
give the cusp cases described in the proposition.}
\label{fig:g1-projections}
\end{figure}
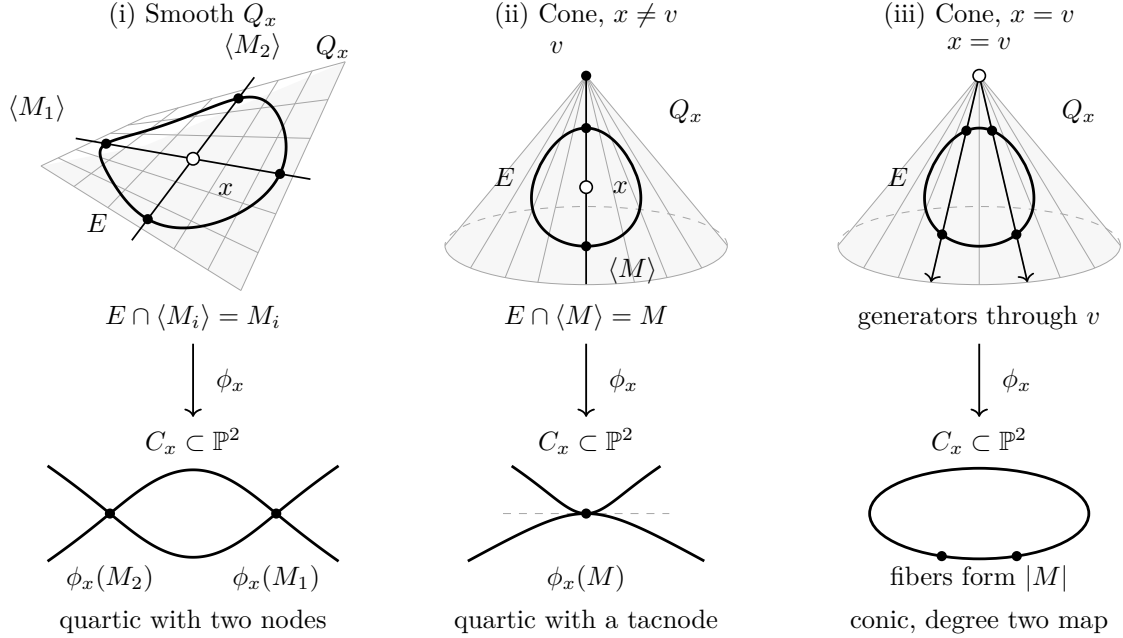

Figure~\ref{fig:g1-projections} illustrates the three cases.

\Needspace{7\baselineskip}
\begin{proof}
We work over \(\overline{\F}_q\).
Since \(E\) is the base locus of \(\Lambda_L\) and \(x\notin E\), a
unique member \(Q_x\) contains \(x\).  Every singular member is a cone,
since a union of two planes or a double plane cannot contain the
irreducible nonplanar curve \(E\).
We first show that the fibers of length at least two correspond exactly
to the lines on \(Q_x\) through \(x\), and that each has length two.
A fiber is the intersection \(E\cap\ell\) for a line \(\ell\)
through \(x\).  If this intersection has length at least two, the
equation of \(Q_x\) restricted to \(\ell\) has at least three zeros,
counting multiplicities and the additional point \(x\notin E\).
A nonzero quadratic on a projective line has only two zeros counted with
multiplicity, so \(\ell\subset Q_x\).
If \(\ell\) lies on \(Q_x\), choose another quadric \(Q'\) in the
pencil.  Since \(E=Q_x\cap Q'\), we have
\(
 E\cap\ell=Q'\cap\ell.
\)
The equation of \(Q'\) restricts to a nonzero quadratic on \(\ell\),
since otherwise \(\ell\) would be contained in \(E\).
Thus this intersection has length two.

In case (i), there are exactly two ruling lines through \(x\).  Their intersections
with \(E\) are degree two divisors \(M_1,M_2\) with disjoint support,
since the lines meet only at \(x\notin E\).
All other nonempty fibers have length one, so \(\phi_x\) is birational
onto a quartic with normalization \(E\) and total delta invariant two.
The length of a normalization fiber equals the multiplicity of its image
point.  Thus the images of \(M_1,M_2\) have multiplicity two, and all
other points have multiplicity one.  A point on a plane curve is smooth
exactly when its multiplicity is one, so these are precisely the two
singular points.  Each has delta invariant one.
A reduced \(M_i\) gives two branches and a doubled \(M_i\) gives one,
so Table~\ref{tab:local-delta-classification} gives a node or a cusp,
respectively.
The plane containing the two ruling lines meets \(E\) in \(M_1+M_2\),
so \(L\simeq\mathcal O_E(M_1+M_2)\).
The planes through \(\langle M_1\rangle\) give the pencil
\(M_1+|M_2|\), and those through \(\langle M_2\rangle\) give
\(M_2+|M_1|\).  Together these span the series of projection from
\(x\).  By \cite[\S13.5, Lemma~13.7]{bergstrom2008cohomology}, this
series gives a degree two map exactly when \([M_1]=[M_2]\).
Our projection is birational, so \([M_1]\neq[M_2]\).

In case (ii), there is exactly one generator through \(x\), namely the line joining
\(x\) and \(v\).  Its intersection with \(E\) is a degree two
fiber \(M\), and all other nonempty fibers have length one.
As in (i), the projection is birational onto a quartic with normalization
\(E\).  The image of \(M\) has multiplicity two, and all other points
have multiplicity one and are smooth.  Hence it is the only singular
point and has delta invariant two.
If \(M=P+Q\) with \(P\neq Q\), its image has two branches
\(B_1,B_2\).  Each has multiplicity one, so both are smooth and have
delta invariant zero.  The formula for the delta invariant of a point
with several branches gives
\[
 2=\delta(B_1)+\delta(B_2)+I(B_1,B_2)=I(B_1,B_2).
\]
Thus the two smooth branches have intersection multiplicity two,
which is the tacnode type \(A_3\) in
Table~\ref{tab:local-delta-classification}.
If \(M=2P\), its image has one branch \(B\), with multiplicity two
and delta invariant two.  In the blowup formula
\[
 \delta(B)=\sum_{i=0}^{N_B}\binom{m_i(B)}2,
\]
the initial multiplicity contributes \(\binom22=1\).
The strict transform after one blowup therefore has delta invariant one.
Its multiplicity must be two, since multiplicity at least three would
already contribute at least three to its delta invariant, while a branch
of multiplicity one is smooth and has delta invariant zero.
One further blowup removes the remaining contribution of one, so the
branch then becomes smooth.  Its multiplicity sequence is therefore
\((2:2:1)\), which is the type \(A_4\) in
Table~\ref{tab:local-delta-classification}.
To obtain the relation for \(L\), project from the vertex \(v\),
which is not on \(E\) because \(E\) is smooth.
The fibers over the base conic form the degree two pencil \(|M|\).
A general plane through \(v\) meets the cone in two generators, whose
intersections with \(E\) are both linearly equivalent to \(M\).
Hence \(L\simeq\mathcal O_E(2M)\).

In case (iii), \(x=v\), and projection maps the cone to its base conic.
Every fiber on \(E\) is its intersection with a generator and has
length two.  The map therefore has degree two, and its fibers form the
complete pencil \(|M|\).  As in (ii), a general plane through \(v\)
gives \(L\simeq\mathcal O_E(2M)\).

Conversely, let \(M_1,M_2\) have disjoint support and
\(L\simeq\mathcal O_E(M_1+M_2)\).
Their spans are distinct lines in the plane whose intersection with
\(E\) is \(M_1+M_2\).  They meet at a point \(x\notin E\), and
both lie on \(Q_x\) by the fiber argument above.
The pencils \(M_1+|M_2|\) and \(M_2+|M_1|\) again span the series
of projection from \(x\).  The degree criterion used in (i) gives a
degree two map exactly when \([M_1]=[M_2]\), which is case (iii).
If \([M_1]\neq[M_2]\), the two distinct lines through \(x\) place
us in case (i).
Finally, suppose \(L\simeq\mathcal O_E(2M)\), and choose a distinct
member \(M'\in|M|\).  Since \(|M|\) is a base point free pencil,
its distinct members \(M,M'\) have disjoint support.  Moreover,
\(M'\sim M\) gives \(L\simeq\mathcal O_E(M+M')\), so the preceding
construction gives the projection center
\(x=\langle M\rangle\cap\langle M'\rangle\notin E\).
Since \([M]=[M']\), the degree criterion shows that projection from
\(x\) has degree two.  Cases (i) and (ii) give birational projections,
so this must be case (iii).  Thus \(Q_x\) is a cone and its vertex is
\(v=x\).  The fibers of projection from \(v\) form the complete
pencil \(|M|\), since \(M\) is one of these fibers.  Consequently, the generators are exactly the
lines \(\langle M'\rangle\) for \(M'\in|M|\), and all these lines
pass through the same vertex \(v\).  Case (ii) now shows that the centers on
\(\langle M\rangle\) outside \(\operatorname{Supp}(M)\cup\{v\}\)
give precisely the projections with singular fiber \(M\).
These constructions are defined over \(\F_q\) when \(M\), or the
unordered pair \(\{M_1,M_2\}\), is preserved by Frobenius.
\end{proof}

The proposition gives the five geometric singularity configurations
\(A_1+A_1\), \(A_1+A_2\), \(A_2+A_2\), \(A_3\), and \(A_4\).
Taking the Frobenius action into account gives the following twelve profiles.
In the two conjugate node rows, \([2^2]\) means that the four geometric
branches form two orbits of size two, and \([4]\) means that they form one
orbit of size four.

\begingroup
\small
\renewcommand{\arraystretch}{1.06}
\setlength{\LTpre}{6pt}
\setlength{\LTpost}{6pt}
\renewcommand{\LTcaptype}{unnumtable}
\begin{longtable}{@{}T{2.0cm}T{6.2cm}T{2.2cm}T{3.5cm}@{}}
\toprule
geometric type & \(\F_q\)-profile & trace shift \(\Delta_\tau\)
& \(\#C_x(\F_q)\) \\
\midrule
\endfirsthead
\toprule
geometric type & \(\F_q\)-profile & trace shift \(\Delta_\tau\)
& \(\#C_x(\F_q)\) \\
\midrule
\endhead
\midrule
\multicolumn{4}{r}{Continued on next page}\\
\endfoot
\bottomrule
\endlastfoot
\(A_1+A_1\) & two rational split nodes & 2 & \(\#E(\F_q)-2\)\\
\(A_1+A_1\) & one split and one nonsplit rational node & 0 & \(\#E(\F_q)\)\\
\(A_1+A_1\) & two rational nonsplit nodes & \(-2\) & \(\#E(\F_q)+2\)\\
\(A_1+A_1\) & conjugate nodes, with branch orbits \([2^2]\) & 0 & \(\#E(\F_q)\)\\
\(A_1+A_1\) & conjugate nodes, with branch orbit \([4]\) & 0 & \(\#E(\F_q)\)\\
\addlinespace
\(A_1+A_2\) & split rational node and rational cusp & 1 & \(\#E(\F_q)-1\)\\
\(A_1+A_2\) & nonsplit rational node and rational cusp & \(-1\) & \(\#E(\F_q)+1\)\\
\addlinespace
\(A_2+A_2\) & two rational cusps & 0 & \(\#E(\F_q)\)\\
\(A_2+A_2\) & conjugate cusps & 0 & \(\#E(\F_q)\)\\
\addlinespace
\(A_3\) & split tacnode & 1 & \(\#E(\F_q)-1\)\\
\(A_3\) & nonsplit tacnode & \(-1\) & \(\#E(\F_q)+1\)\\
\addlinespace
\(A_4\) & rational ramphoid cusp & 0 & \(\#E(\F_q)\)\\
\end{longtable}
\endgroup

\subsubsection{Counting by trace shift}

Choose an origin \(O\in E(\F_q)\), and identify \(E(\F_q)\) with
\(\Pic^2(E)(\F_q)\) by \(P\mapsto[P+O]\).
We write \(2E(\F_q)\) for the image of multiplication by two and
\(E[2](\F_q)\) for its kernel.  Under this identification,
\(a\in2E(\F_q)\) means that \(a=[2P]\) for some \(P\in E(\F_q)\).
For each \(a\in\Pic^2(E)(\F_q)\), the pencil \(|a|\) has
\(q+1\) rational members.
A rational member is either a pair of
distinct rational points, a doubled rational point, or a pair of conjugate
points in \(E(\F_{q^2})\setminus E(\F_q)\).  These are called split
reduced, doubled, and nonsplit divisors, respectively.

\Needspace{12\baselineskip}
\begin{prop}\label{prop:g1-degree-two-divisors}
For each \(a\in\Pic^2(E)(\F_q)\), the numbers of split reduced,
doubled, and nonsplit divisors in \(|a|\) are as follows.
\begin{center}
\normalfont
\renewcommand{\arraystretch}{1.6}
\begin{tabular}{@{}lcc@{}}
\toprule
divisor type & \(a\in2E(\F_q)\) & \(a\notin2E(\F_q)\) \\
\midrule
split reduced
& \(\displaystyle\frac{\#E(\F_q)-\#E[2](\F_q)}2\)
& \(\displaystyle\frac{\#E(\F_q)}2\)\\
doubled & \(\#E[2](\F_q)\) & \(0\)\\
nonsplit
& \(\displaystyle\frac{2q+2-\#E(\F_q)-\#E[2](\F_q)}2\)
& \(\displaystyle\frac{2q+2-\#E(\F_q)}2\)\\
\bottomrule
\end{tabular}
\end{center}
Exactly \(\#E(\F_q)/\#E[2](\F_q)\) of the rational divisor classes
lie in \(2E(\F_q)\).
Across all classes, the total numbers of split reduced and doubled
divisors are \(\binom{\#E(\F_q)}2\) and \(\#E(\F_q)\), respectively.
The total number of nonsplit divisors is
\[
 \frac{\#E(\F_{q^2})-\#E(\F_q)}2
 =\frac{\#E(\F_q)(2q+1-\#E(\F_q))}2.
\]
\end{prop}

\begin{proof}
For each \(P\in E(\F_q)\) there is a unique \(Q\in E(\F_q)\)
with \([P+Q]=a\), so there are \(\#E(\F_q)\) ordered pairs.
Those with \(P=Q\) correspond to the doubled divisors.
Their number is \(\#E[2](\F_q)\) when \(a\in2E(\F_q)\), and zero
otherwise.  Removing these pairs and dividing by two gives the split
reduced count.  Subtracting the split reduced and doubled counts from
\(q+1\) gives the nonsplit count.

Since multiplication by two has kernel \(E[2](\F_q)\), exactly
\(\#E(\F_q)/\#E[2](\F_q)\) of the rational divisor classes lie in
\(2E(\F_q)\).

Across all classes, counting unordered pairs of distinct rational
points gives \(\binom{\#E(\F_q)}2\) split reduced divisors, and
counting rational points gives \(\#E(\F_q)\) doubled divisors.
Each nonsplit divisor is an unordered pair of conjugate points in
\(E(\F_{q^2})\setminus E(\F_q)\), so their total number is
\((\#E(\F_{q^2})-\#E(\F_q))/2\).
These divisors also correspond to the closed points of degree two on \(E\).
The second expression for their number follows from the relation
\(
 \#E(\F_{q^2})=(q+1)^2-(q+1-\#E(\F_q))^2.
\)
\end{proof}

The twelve genus one profiles listed above have trace shifts in
\(\{-2,-1,0,1,2\}\).  We group the plane series on the fixed curve \(E\)
by their trace shift.  For each \(d\) in this set, define
\[
 B_d(E)=\#\{\gamma\in\mathcal G^2_4(E):\Delta(C_\gamma)=d\}.
\]

\begin{prop}\label{prop:g1-fixed-counts}
The count for trace shift two is
\[
 B_2(E)=\frac{\#E(\F_q)}8\bigl(
 \#E(\F_q)^3-7\#E(\F_q)^2+15\#E(\F_q)-\#E[2](\F_q)-8\bigr).
\]
The total for trace shift one is
\[
 B_1(E)=\frac{\#E(\F_q)}2\bigl(
 q\#E(\F_q)+\#E(\F_q)^2-q-6\#E(\F_q)+\#E[2](\F_q)+4\bigr).
\]
For trace shift minus one, we obtain
\[
 B_{-1}(E)=\frac{\#E(\F_q)}2\bigl(
 2q^2+q\#E(\F_q)-\#E(\F_q)^2-q+2\#E(\F_q)
 +\#E[2](\F_q)-2\bigr),
\]
and trace shift minus two gives
\[
\begin{aligned}
 B_{-2}(E)=\frac{\#E(\F_q)}8\bigl(&
 4q^2\#E(\F_q)-4q\#E(\F_q)^2+\#E(\F_q)^3-4q^2\\
 &+8q\#E(\F_q)-3\#E(\F_q)^2-4q+3\#E(\F_q)-\#E[2](\F_q)
 \bigr).
\end{aligned}
\]
The count for trace shift zero is
\[
\begin{aligned}
 B_0(E)=\frac{\#E(\F_q)}4\bigl(&
 4q^3-2q^2\#E(\F_q)+2q\#E(\F_q)^2-\#E(\F_q)^3+2q^2\\
 &-8q\#E(\F_q)+5\#E(\F_q)^2+10q-5\#E(\F_q)-3\#E[2](\F_q)
 \bigr).
\end{aligned}
\]
Their sum is
\[
 \sum_{d=-2}^{2}B_d(E)
 =\#E(\F_q)(q^3+q^2+q-\#E(\F_q)).
\]
\end{prop}
\begin{proof}
First count the series giving quartics with one geometric singular point.
By Proposition~\ref{prop:g1-projection-geometry}(ii), we choose a rational
effective divisor \(M\) of degree two and then a rational center on
\(\langle M\rangle\), excluding the cone vertex and
\(\operatorname{Supp}(M)\).
There are \(q-2\), \(q-1\), and \(q\) such centers when \(M\) is split
reduced, doubled, and nonsplit, respectively.
Proposition~\ref{prop:g1-degree-two-divisors} therefore gives
\(\binom{\#E(\F_q)}2(q-2)\) split tacnodes,
\(\#E(\F_q)(q-1)\) ramphoid cusps, and
\(q(\#E(\F_{q^2})-\#E(\F_q))/2\) nonsplit tacnodes.
Each series is counted once, since it determines \(M\) and the projection
center.

For two geometric singular points, use the divisor pairs of
Proposition~\ref{prop:g1-projection-geometry}(i).
Frobenius must preserve the pair.  If it exchanges the two divisors,
the trace shift is zero.  For two rational singular points, the profiles
with nonzero trace shift have the following counts.  In each case,
subtract the entry for equal divisor classes from the entry for disjoint
support.
\begingroup
\small
\renewcommand{\arraystretch}{1.15}
\setlength{\LTpre}{6pt}
\setlength{\LTpost}{6pt}
\renewcommand{\LTcaptype}{unnumtable}
\begin{longtable}{@{}T{3.3cm}T{2.8cm}T{7.4cm}@{}}
\toprule
\(\F_q\)-profile & \mbox{divisor pair} \mbox{condition} & number\\
\midrule
\endfirsthead
\toprule
\(\F_q\)-profile & \mbox{divisor pair} \mbox{condition} & number\\
\midrule
\endhead
\midrule
\multicolumn{3}{r}{Continued on next page}\\
\endfoot
\bottomrule
\endlastfoot
two split nodes
& disjoint & \(3\binom{\#E(\F_q)}4\)\\*
two split nodes
& equal divisor class
& \(\displaystyle\frac{\#E(\F_q)
 (\#E(\F_q)^2-4\#E(\F_q)+\#E[2](\F_q)+2)}8\)\\[1mm]
\addlinespace
split node and cusp
& disjoint & \(\#E(\F_q)\binom{\#E(\F_q)-1}2\)\\*
split node and cusp
& equal divisor class
& \(\displaystyle\frac{\#E(\F_q)(\#E(\F_q)-\#E[2](\F_q))}2\)\\[1mm]
\addlinespace
nonsplit node and cusp
& disjoint
& \(\displaystyle\frac{\#E(\F_q)(\#E(\F_{q^2})-\#E(\F_q))}2\)\\*
nonsplit node and cusp
& equal divisor class
& \(\displaystyle\frac{\#E(\F_q)(2q+2-\#E(\F_q)-\#E[2](\F_q))}2\)\\[1mm]
\addlinespace
two nonsplit nodes
& disjoint
& \(\displaystyle\binom{(\#E(\F_{q^2})-\#E(\F_q))/2}{2}\)\\*
two nonsplit nodes
& equal divisor class
& \(\displaystyle\frac{\#E(\F_q)
 ((2q-\#E(\F_q))^2+\#E[2](\F_q)-2)}8\)\\
\end{longtable}
\endgroup
The counts for disjoint support are obtained by choosing the points in
the two divisors.  For example, there are three ways to divide four
rational points into two unordered pairs, giving
\(3\binom{\#E(\F_q)}4\).
For the subtraction terms, apply Proposition~\ref{prop:g1-degree-two-divisors}
in each divisor class.  Choose two distinct divisors of the same type
for two nodes of the same kind, or one divisor of each type for a node
and a cusp.  Distinct members of \(|a|\) have disjoint support, so no
further restriction is needed.  For example, the term for two split
nodes is
\[
\begin{aligned}
 &\frac{\#E(\F_q)}{\#E[2](\F_q)}
   \left[\binom{(\#E(\F_q)-\#E[2](\F_q))/2}{2}
   +(\#E[2](\F_q)-1)\binom{\#E(\F_q)/2}{2}\right]\\
 &=\frac{\#E(\F_q)
   (\#E(\F_q)^2-4\#E(\F_q)+\#E[2](\F_q)+2)}8.
\end{aligned}
\]
The other three subtraction terms follow in the same way.

The counts for two split nodes and two nonsplit nodes give \(B_2(E)\)
and \(B_{-2}(E)\), respectively.
For \(B_1(E)\), add the split tacnode count to the count for a split
node and a cusp.  For \(B_{-1}(E)\), add the nonsplit tacnode count to
the count for a nonsplit node and a cusp.

For the total count, there are \(\#E(\F_q)\) rational line bundles
\(L\) of degree four.  For each \(L\), the centers outside \(E\) give
\(q^3+q^2+q+1-\#E(\F_q)\) base point free series.
Each nonbirational series comes from a unique
\(A\in\Pic^2(E)(\F_q)\), with \(L=A^{\otimes2}\) and subspace
\(\operatorname{im}(\Sym^2H^0(E,A)\to H^0(E,A^{\otimes2}))\).
There are therefore \(\#E(\F_q)\) such conic series, so the number of
birational series is
\[
 \#E(\F_q)(q^3+q^2+q+1-\#E(\F_q))-\#E(\F_q)
 =\#E(\F_q)(q^3+q^2+q-\#E(\F_q)).
\]
Subtracting the four nonzero trace shift counts gives \(B_0(E)\).
The divisor counts and their polynomial simplifications are recorded in
the accompanying code.
\end{proof}

Combining Propositions~\ref{prop:g1-fixed-counts} and
\ref{prop:irr-strata-contrib} gives the contribution of quartics with
genus one normalization:
\begin{equation}\label{eq:g1-weight-enumerator}
   (q-1)\#\PGL_3(\F_q)
   \sum_{\substack{[E]\;\text{smooth genus one}/\F_q}}
   \frac{1}{\#\Aut_{\F_q}(E)}
   \sum_{d=-2}^{2}
   B_d(E)
   X^{\#E(\F_q)-d}
   Y^{q^2+q+1-\#E(\F_q)+d}.
\end{equation}
Here the sum is over unpointed smooth genus one curves and
\(\Aut_{\F_q}(E)\) is the full automorphism group of the underlying curve,
including rational translations.  The choice of origin does not affect
\(\#E[2](\F_q)\), so it does not affect any \(B_d(E)\).

\subsection{Geometric genus zero}
\label{chap:genus-zero}

Let \(C_0\subset\PP^2\) be a geometrically integral plane quartic whose
normalization has genus zero.  Over a finite field every smooth projective
geometrically connected curve of genus zero is isomorphic to \(\PP^1\).
Let \(\nu\colon\PP^1\longrightarrow C_0\) be the normalization map, and
let \(j\colon C_0\hookrightarrow\PP^2\) be the inclusion.  Write
\(\phi=j\circ\nu\).  Then
\(\phi^*\mathcal O_{\PP^2}(1)=\nu^*\mathcal O_{C_0}(1)
\simeq\mathcal O_{\PP^1}(4)\).
The pullbacks of the coordinate linear forms \(X_0,X_1,X_2\) on
\(\PP^2\) span a three dimensional base point free subspace
\(V\subset H^0(\PP^1,\mathcal O_{\PP^1}(4))\),
and \(\phi\) is the morphism defined by \(V\).  Conversely, suppose
that the morphism defined by such a subspace is birational onto its image.  It then gives a geometrically
integral rational plane quartic, unique up to projective coordinates.

For a Galois refined profile \(\tau\), the number \(m_\tau(\PP^1)\)
counts the subspaces \(V\) described above.  Define
\[
  \overline m_\tau
  =\frac{m_\tau(\PP^1)}{\#\PGL_2(\F_q)}.
\]
By Proposition~\ref{prop:irr-strata-contrib}, the number of embedded
quartics of profile \(\tau\) is
\(\#\PGL_3(\F_q)\overline m_\tau\).
We call \(\overline m_\tau\) the normalized count.  It is a weighted
count of the orbits under changes of parameter on \(\PP^1\), with
each orbit weighted by the reciprocal of its stabilizer order.
We compute these counts and group them by trace shift to obtain the
genus zero contribution.

The total delta invariant is three, so a genus zero quartic has one,
two, or three geometric singular points.  For curves with two or three
singular points, we use the counts in
\cite[\S\S12.3 and~13.4]{bergstrom2008cohomology} and recall the
divisor constructions that produce them.  For curves with one singular
point, we count subspaces directly and use local calculations to
determine the singularity type.

To construct quartics with three geometric singular points, let
\(\{M_1,M_2,M_3\}\) be an unordered triple of effective degree two
divisors on \(\PP^1_{\overline{\F}_q}\), preserved by Frobenius and with
pairwise disjoint support.  Choose \(s_i\in H^0(\PP^1,\mathcal O_{\PP^1}(2))\) with
\((s_i)_0=M_i\), and set
\(V
   =\langle s_1s_2,s_1s_3,s_2s_3\rangle
   \subset H^0(\PP^1,\mathcal O_{\PP^1}(4))\).
Since \((s_is_j)_0=M_i+M_j\), the projective linear system \(\PP(V)\)
is the one spanned by \(M_1+M_2,M_1+M_3,M_2+M_3\) in
\cite[\S12.3]{bergstrom2008cohomology}.
The three products have no common zero.
Over \(\overline{\F}_q\), the displayed basis identifies
\(\PP(V^\vee)\) with \(\PP^2\), and \(\phi\) sends
\(M_1,M_2,M_3\) to the three coordinate points.
Since \(\phi^*\mathcal O_{\PP^2}(1)\simeq\mathcal O_{\PP^1}(4)\),
\[
  4=\deg(\phi\colon\PP^1\to C_0)\deg(C_0).
\]
The image contains three noncollinear points, so it is not a line.
Thus \(\phi\) is either birational onto a quartic or has degree two
onto a conic.

The linear system spanned by \(M_1,M_2,M_3\) in
\cite[Lemma~12.8]{bergstrom2008cohomology} is
\(\PP(\langle s_1,s_2,s_3\rangle)\).  Its projective dimension is
\(\dim\langle s_1,s_2,s_3\rangle-1\), so it is a pencil exactly when
\(\dim\langle s_1,s_2,s_3\rangle=2\).  It is base point free because
the \(M_i\) have disjoint support.  Thus Bergstr\"om's
Lemma~12.8 \cite{bergstrom2008cohomology} says that
\(\phi\) has degree two exactly in this case.
Otherwise \(\phi\) is birational, and
\cite[Lemma~12.7]{bergstrom2008cohomology} shows that \(C_0\) has
three singularities.  Their delta invariants sum to three, so each
has delta invariant one.  After excluding the conic maps,
the normalized counts for all 22 Frobenius profiles are given in
\cite[Lemma~12.13 and Table~6]{bergstrom2008cohomology}.
We use the odd characteristic column of the table.

For a quartic with two geometric singular points, their delta invariants
are two and one.  Let \(M_1\) and \(M_2\) be their respective scheme
theoretic inverse images under the normalization map \(\nu\).  Since their
delta invariants differ, both singular points are rational, and
\(M_1,M_2\in(\Sym^2\PP^1)(\F_q)\) have disjoint support.  The number
\(\#(\operatorname{Supp}(M_1)\cap\PP^1(\F_q))\) is one for a ramphoid cusp,
two for a split tacnode, and zero for a nonsplit tacnode.  For fixed
\(M_1,M_2\), the construction in
\cite[\S13.4, Lemma~13.4]{bergstrom2008cohomology} gives exactly
\(q-\#(\operatorname{Supp}(M_1)\cap\PP^1(\F_q))\) base point free
linear series birational onto quartics with these singular fibers.
Counting all such pairs and dividing by \(\#\PGL_2(\F_q)\) gives
the normalized counts in \cite[Lemma~13.5]{bergstrom2008cohomology}.
Translating the singularity notation and grouping
by trace shift gives the following normalized counts for nonzero shifts.
These trace shift counts are calculated in the accompanying code.
\begin{center}
\captionsetup{hypcap=false}
\captionof{table}{Normalized counts for quartics with two or three
geometric singular points and nonzero trace shift.}
\label{tab:g0-two-three-singular-counts}
\renewcommand{\arraystretch}{1.22}
\begin{tabular}{@{}cT{5.1cm}T{5.1cm}@{}}
\toprule
trace shift \(\Delta_\tau\) & three singular points & two singular points\\
\midrule
 \(3\)&\(\dfrac{(q-3)(q^2-7q+11)}{48}\)&0\\[2mm]
 \(2\)&\(\dfrac{(q-3)^2}{8}\)&\(\dfrac{(q-2)^2}{4}\)\\[2mm]
 \(1\)&\(\dfrac{5q^3-8q^2+12q-13}{16}\)&\(\dfrac{2q-3}{2}\)\\[2mm]
\(-1\)&\(\dfrac{5q^3-14q^2+8q+3}{16}\)&\(\dfrac{2q-1}{2}\)\\[2mm]
\(-2\)&\(\dfrac{(q-1)^2}{8}\)&\(\dfrac{q(q-2)}{4}\)\\[2mm]
\(-3\)&\(\dfrac{(q-1)(q^2-3q+1)}{48}\)&0\\
\bottomrule
\end{tabular}
\end{center}

If there is only one geometric singular point, it is \(\F_q\)-rational
and has delta invariant three.  Its geometric type is one of
\(D_4,D_5,E_6,A_5,A_6\), as described in
Table~\ref{tab:local-delta-classification}.

For \(D_4\), the superscript records the lengths of the Frobenius orbits
on its three branches.  We write \(D_4^{1+1+1}\) when all three branches
are fixed, \(D_4^{1+2}\) when one is fixed and the other two are exchanged,
and \(D_4^3\) when Frobenius cyclically permutes all three.
For \(A_5\), the superscripts \(\mathrm{sp}\) and \(\mathrm{ns}\) mean
split and nonsplit, respectively.  In \(A_5^{\mathrm{sp}}\), both branches
are fixed individually, whereas in \(A_5^{\mathrm{ns}}\) they are exchanged.
The two branches of \(D_5\) have different multiplicities, so Frobenius
fixes each one.  The types \(E_6\) and \(A_6\) each have a single branch,
which is also fixed.  These give eight profiles.  Since the singular
point is rational, \(\Delta_\tau\) is the number of branches fixed by
Frobenius minus one.

The two constructions below use the following unit count.

\begin{lemm}\label{lem:divisor-unit-count}
Let \(M=\sum_Q e_QQ\) be a nonzero effective divisor on \(\PP^1\)
defined over \(\F_q\), written as a sum of distinct closed points
\(Q\) with multiplicities \(e_Q\).  Then
\[
 \#H^0(\mathcal O_M)^\times
 =\prod_Q q^{(e_Q-1)\deg Q}(q^{\deg Q}-1).
\]
\end{lemm}

\begin{proof}
The local factor at \(Q\) has \(q^{e_Q\deg Q}\) elements and
residue field \(\F_{q^{\deg Q}}\).  Its maximal ideal has
\(q^{(e_Q-1)\deg Q}\) elements, so it has
\(q^{(e_Q-1)\deg Q}(q^{\deg Q}-1)\) units.
The algebra \(H^0(\mathcal O_M)\) is the product of these local
factors, giving the formula.
\end{proof}

For a series with one singular point \(P\), let
\(M=\nu^{-1}(P)\) be its scheme theoretic fiber.
The pullbacks of linear forms vanishing at \(P\) form a two dimensional
subspace \(W\subset V\) whose common zero divisor \(M\) records
the initial branch multiplicities.
Table~\ref{tab:local-delta-classification} gives \(\deg M=3\) for
\(D_4,D_5,E_6\), and \(\deg M=2\) for \(A_5,A_6\).
Choose a section \(s_M\) of \(\mathcal O_{\PP^1}(\deg M)\) with
zero divisor \(M\).  Since \(P\) is rational and \(\nu\) is an
isomorphism away from \(P\), the trace shift is
\(\Delta_\tau=\#M(\F_q)-1\), where \(\#M(\F_q)\) counts rational
support points without multiplicity.

\Needspace{12\baselineskip}
\begin{prop}\label{prop:g0-triple-point-counts}
For profiles with a singularity of type \(D_4,D_5,E_6\), the
normalized counts and trace shifts are
\begin{center}
\renewcommand{\arraystretch}{1.2}
\begin{tabular}{@{}T{4.8cm}cc@{}}
\toprule
\(\F_q\)-profile & trace shift \(\Delta_\tau\) & normalized count \(\overline m_\tau\)\\
\midrule
\(D_4^{1+1+1}\)&2&\(\dfrac{(q-1)^2}{6}\)\\[1mm]
\(D_4^{1+2}\)&0&\(\dfrac{q^2-1}{2}\)\\[1mm]
\(D_4^{3}\)&\(-1\)&\(\dfrac{q^2+q+1}{3}\)\\[1mm]
\(D_5\)&1&\(q\)\\[1mm]
\(E_6\)&0&\(\dfrac{q}{q-1}\)\\
\bottomrule
\end{tabular}
\end{center}
\end{prop}

\begin{proof}
Here \(\deg M=3\), so dimension two forces
\(W=H^0(\mathcal O_{\PP^1}(4)(-M))
=s_MH^0(\mathcal O_{\PP^1}(1))\).
An extension \(V=W+\langle f\rangle\) is base point free exactly
when \(f\) is nonzero on \(\operatorname{Supp}(M)\).
Restriction to \(M\) and a trivialization of \(\mathcal O_M(4)\)
identify the extensions with
\(H^0(\mathcal O_M)^\times/\F_q^\times\): adding an element of
\(W\) or rescaling \(f\) does not change \(V\).
Conversely, every such extension for a rational effective divisor
\(M\) of degree three is birational: the ratio of two generators
of \(W\) is a degree one coordinate.  The image is a rational
quartic with fiber \(M\) over \(P\).  Thus \(m(P)=3\), and
Table~\ref{tab:local-delta-classification} gives \(D_4,D_5,E_6\)
according as \(M\) has three, two, or one geometric support points.
Each has \(\delta_P=3\), so \(P\) is the only singular point.
Its fiber recovers \(M\), and hence \(W\), so each series is counted once.

For the three \(D_4\)-profiles in the stated order, there are
\(\binom{q+1}{3}\), \((q+1)(q^2-q)/2\), and \((q^3-q)/3\)
choices of \(M\).  By Lemma~\ref{lem:divisor-unit-count}, dividing
the unit counts by \(q-1\) gives respectively \((q-1)^2\),
\(q^2-1\), and \(q^2+q+1\) series per fiber.
For \(D_5\), there are \(q(q+1)\) fibers \(M=2p+r\) with
distinct rational points \(p,r\), each giving \(q(q-1)\) series.
For \(E_6\), there are \(q+1\) fibers \(M=3p\) with rational
\(p\), each giving \(q^2\) series.  Multiplying these fiber and
series counts and dividing by \(\#\PGL_2(\F_q)\) gives the
normalized counts.  The rational support of \(M\) gives the stated
trace shifts.
\end{proof}

\Needspace{9\baselineskip}
\begin{prop}\label{prop:g0-a5-a6-counts}
For profiles with a singularity of type \(A_5,A_6\), the normalized
counts and trace shifts are
\begin{center}
\renewcommand{\arraystretch}{1.2}
\begin{tabular}{@{}T{4.8cm}cc@{}}
\toprule
\(\F_q\)-profile & trace shift \(\Delta_\tau\) & normalized count \(\overline m_\tau\)\\
\midrule
\(A_5^{\mathrm{sp}}\)&1&\(\dfrac{q-1}{2}\)\\[1mm]
\(A_5^{\mathrm{ns}}\)&\(-1\)&\(\dfrac{q-1}{2}\)\\[1mm]
\(A_6\)&0&1\\
\bottomrule
\end{tabular}
\end{center}
\end{prop}

\begin{proof}
Here \(\deg M=2\).  We first determine the necessary form of the series.
The tangent line pulls back to \(2M\).  Indeed, for \(A_5\) write
\(M=p_1+p_2\).  The contact is at least two with each smooth branch,
so the pullback contains \(2p_1+2p_2\).
For \(A_6\), contact three would give orders \((2,3)\)
and a smooth branch after one blowup, contrary to the sequence
\((2:2:2:1)\).  Its contact is therefore at least four.  In both cases
the pullback has degree four, forcing equality.
Removing the common zero divisor from \(W\) gives
\(W=s_MU\), where \(U\subset H^0(\mathcal O_{\PP^1}(2))\) is a base
point free pencil.  The tangent line passes through \(P\), so its
pullback is a nonzero scalar multiple of \(s_M^2\) in \(W\).
Hence \(s_M\in U\).  Choose \(h\) with \(U=\langle s_M,h\rangle\)
and write \(V=\langle s_M^2,s_Mh,f\rangle\).

In the chart where \(f\ne0\), put \(u=s_Mh/f\) and \(v=s_M^2/f\).
Then \(v/u^2=f/h^2\) and \(v-u^2=s_M^2(f-h^2)/f^2\).
If \(M\) is reduced, \(u\) is a local parameter on each branch.
Their intersection multiplicity is three, so their quadratic
coefficients, the values of \(f/h^2\) at the two points, agree.
The common nonzero value is Frobenius fixed, and rescaling \(f\) over
\(\F_q\) makes it one.  If \(M=2p\), similarly normalize
\(f(p)=h(p)^2\).  Were \(f-h^2\) to vanish only once at \(p\),
the coordinates \(u,v-u^2\) would have orders \((2,5)\), hence
multiplicity sequence \((2:2:1)\), contradicting the \(A_6\) sequence.
Thus in both cases \(f-h^2\in s_MH^0(\mathcal O_{\PP^1}(2))\).
Choose \(g\in H^0(\mathcal O_{\PP^1}(2))\setminus U\).
When \(M=2p\), subtract a multiple of \(h\) to arrange \(g(p)=0\).
This zero is simple, since otherwise \(g\) would be a multiple of
\(s_M\).  Reducing \((f-h^2)/s_M\) modulo \(U\) shows that every
relevant series has the form
\[
 V=\langle s_M^2,s_Mh,h^2+\lambda s_Mg\rangle,
 \qquad \lambda\in\F_q.
\]

Conversely, fix \(M,U,h,g\) as above.
The value \(\lambda=0\) gives \(\Sym^2U\), a degree two map onto a
conic, and is excluded.  For \(\lambda\ne0\), put
\(f=h^2+\lambda s_Mg\).  The series is base point free because
\(s_M^2\) is nonzero away from \(M\) and \(f=h^2\) is nonzero on
its support.  The identity \(v-u^2=\lambda s_M^3g/f^2\) gives
\(s_M/h=v/u\) and \(g/h=(v-u^2)/(\lambda uv)\).
These ratios recover the quadratic embedding defined by the basis
\(s_M,h,g\) of \(H^0(\mathcal O_{\PP^1}(2))\).
Thus the map is birational onto a quartic.
The fiber over \([0:0:1]\) is exactly \(M\), since its local ideal
\((s_Mh/f,s_M^2/f)\) is \((s_M)\).

If \(M\) is reduced, the two branches are smooth and
\(v-u^2\) vanishes to order at least three on each.
Their intersection multiplicity is therefore at least three, and
\eqref{eq:local-delta-branches} and the total delta invariant force it
to be exactly three.  This is the \(A_5\) condition in
Table~\ref{tab:local-delta-classification}.
If \(M=2p\), the coordinates \(u,v-u^2\) have orders \((2,7)\),
since \(g\) vanishes simply at \(p\).  Successive blowups give
\((2,5)\), \((2,3)\), and \((2,1)\), hence the \(A_6\) sequence
\((2:2:2:1)\).  In both cases the delta invariant is three, so there
are no other singularities.  Every nonzero parameter therefore gives
a valid series of the required type.

These parameters do not overcount.  For fixed \(h,g\), equality of the
subspaces for \(\lambda,\mu\) forces the relative scalar of the third
generator to be one by restriction to \(M\).  It follows that
\((\lambda-\mu)g\in U\), and \(g\notin U\) gives \(\lambda=\mu\).
The series also recovers its unique singular fiber \(M\), the pencil
\(W\) of pullbacks through its image, and hence \(U\).
For fixed \(M\), restriction to \(M\) identifies the base point free
pencils containing \(s_M\) with
\(H^0(\mathcal O_M)^\times/\F_q^\times\), after trivializing
\(\mathcal O_M(2)\).  Each pencil gives \(q-1\) distinct valid
extensions, so the number of series per fiber is
\(\#H^0(\mathcal O_M)^\times\).
There are \(\binom{q+1}{2}\) split reduced fibers,
\((q^2-q)/2\) nonsplit fibers, and \(q+1\) doubled fibers.
Lemma~\ref{lem:divisor-unit-count} gives respectively \((q-1)^2\),
\(q^2-1\), and \(q(q-1)\) series per fiber.
Multiplying and dividing by \(\#\PGL_2(\F_q)\) gives the stated
normalized counts.  The rational support of \(M\) gives the trace shifts.
\end{proof}

\subsubsection{Counting by trace shift}
\label{subsec:g0-trace-shift-counts}

We group the linear series on the fixed normalization \(\PP^1\) by
their trace shift.  For \(-3\le d\le3\), define
\[
 C_d(q)=\#\{\gamma\in\mathcal G^2_4(\PP^1):\Delta(C_\gamma)=d\}.
\]

\begin{prop}\label{prop:g0-trace-shift-counts}
The trace shift counts are
\begin{align*}
 C_3(q)&=\frac{q(q^2-1)(q-3)(q^2-7q+11)}{48},\\
 C_2(q)&=\frac{q(q^2-1)(13q^2-50q+55)}{24},\\
 C_1(q)&=\frac{q(q^2-1)(5q^3-8q^2+52q-45)}{16},\\
 C_0(q)&=\frac{q(q+1)(4q^4+13q^3-35q^2+49q-19)}{12},\\
 C_{-1}(q)&=\frac{q(q^2-1)(15q^3-26q^2+112q-23)}{48},\\
 C_{-2}(q)&=\frac{q(q^2-1)(3q^2-6q+1)}{8},\\
 C_{-3}(q)&=\frac{q(q^2-1)(q-1)(q^2-3q+1)}{48},\\
 \sum_{d=-3}^3 C_d(q)&=q^2(q+1)(q^3+q-1).
\end{align*}
\end{prop}

\begin{proof}
We use the normalized counts in Table~\ref{tab:g0-two-three-singular-counts},
obtained from \cite[Table~6 and Lemma~13.5]{bergstrom2008cohomology},
and the counts for one singular point in
Propositions~\ref{prop:g0-triple-point-counts} and~\ref{prop:g0-a5-a6-counts}.
For shifts \(3\) and \(-3\), only quartics with three geometric singular
points contribute.  Thus
\[
 \frac{C_3(q)}{\#\PGL_2(\F_q)}
 =\frac{(q-3)(q^2-7q+11)}{48},
 \qquad
 \frac{C_{-3}(q)}{\#\PGL_2(\F_q)}
 =\frac{(q-1)(q^2-3q+1)}{48}.
\]
For shift \(2\), we add the counts for three and two singular points
in Table~\ref{tab:g0-two-three-singular-counts} and the count for
\(D_4^{1+1+1}\) in Proposition~\ref{prop:g0-triple-point-counts}.
For shift \(-2\), only the two table entries contribute.  This gives
\[
 \frac{C_2(q)}{\#\PGL_2(\F_q)}
 =\frac{13q^2-50q+55}{24},
 \qquad
 \frac{C_{-2}(q)}{\#\PGL_2(\F_q)}
 =\frac{3q^2-6q+1}{8}.
\]
For each of the shifts \(1\) and \(-1\), we again add the counts for
three and two singular points in Table~\ref{tab:g0-two-three-singular-counts}.
The additional profiles with one singular point are \(D_5\) and
\(A_5^{\mathrm{sp}}\) for shift \(1\), and \(D_4^3\) and
\(A_5^{\mathrm{ns}}\) for shift \(-1\).
Propositions~\ref{prop:g0-triple-point-counts} and~\ref{prop:g0-a5-a6-counts}
therefore give
\begin{align*}
 \frac{C_1(q)}{\#\PGL_2(\F_q)}
 &=\frac{5q^3-8q^2+52q-45}{16},\\
 \frac{C_{-1}(q)}{\#\PGL_2(\F_q)}
 &=\frac{15q^3-26q^2+112q-23}{48}.
\end{align*}
Multiplying by \(\#\PGL_2(\F_q)=q(q^2-1)\) gives
the six nonzero trace shift formulas.

We compute the total by counting the elements of \(\mathcal G^2_4(\PP^1)\).  A
three dimensional subspace of \(H^0(\PP^1,\mathcal O_{\PP^1}(4))\) can have a
base divisor of degree at most two.  There are
\(
   \#\Gr(3,5)(\F_q)
   =q^6+q^5+2q^4+2q^3+2q^2+q+1
\)
such subspaces in total.  Those with base divisor of degree exactly two
correspond to \((\Sym^2\PP^1)(\F_q)\), so there are \(q^2+q+1\) of them.
The number with base divisor of degree exactly one is
\(
   (q+1)\bigl(\#\Gr(3,4)(\F_q)-(q+1)\bigr)
   =(q+1)(q^3+q^2).
\)
It follows that the number of base point free subspaces is
\(
   q^4(q^2+q+1).
\)
A nonbirational map has degree two onto a conic.  Its linear series has the
unique form
\(\Sym^2U\), where
\(U\subset H^0(\PP^1,\mathcal O_{\PP^1}(2))\) is a base point free pencil.  There
are
\(
   \#\Gr(2,3)(\F_q)-(q+1)=q^2
\)
such pencils.  Therefore the number of birational plane series is
\(q^6+q^5+q^4-q^2=q^2(q+1)(q^3+q-1)\).
The remaining count is
\[
\begin{aligned}
 C_0(q)&=q^2(q+1)(q^3+q-1)\\
 &\quad-C_3(q)-C_2(q)-C_1(q)-C_{-1}(q)-C_{-2}(q)-C_{-3}(q).
\end{aligned}
\]
Substituting the six nonzero counts gives the stated formula.
The profile sums and this subtraction are carried out in the accompanying
code.
\end{proof}

Combining Propositions~\ref{prop:g0-trace-shift-counts} and
\ref{prop:irr-strata-contrib} gives the contribution
\(W^{\Irr,g=0}_{\cC_{2,4}}\) of quartics with genus zero normalization
to the weight enumerator:
\begin{equation}\label{eq:g0-seven-term-weight}
   W^{\Irr,g=0}_{\cC_{2,4}}(X,Y)
   =\frac{(q-1)\#\PGL_3(\F_q)}{\#\PGL_2(\F_q)}
    \sum_{d=-3}^{3}
    C_d(q)X^{q+1-d}Y^{q^2+d}.
\end{equation}

\section{The weight enumerator}
\label{sec:weight-enumerator}

The weight enumerator decomposes as
\[
 W_{\cC_{2,4}}(X,Y)
 =X^{q^2+q+1}+W_{\mathrm{red}}(X,Y)
  +W_{\mathrm{irr,sing}}(X,Y)+W_{\mathrm{sm}}(X,Y).
\]
We expand these terms using the counts from Sections~2 and~3, then
express the remaining genus one and genus two sums in arithmetic terms.
All isomorphisms and automorphisms below are defined over \(\F_q\).
For \(-2\le d\le2\), collect the genus one counts \(B_d(E)\) from
Proposition~\ref{prop:g1-fixed-counts} at each trace as
\[
   \mathsf B_{q,d}(t)
   =
   \sum_{\substack{[E]\;\mathrm{smooth\ genus\ one}/\F_q\\t_E=t}}
   \frac{B_d(E)}{\#\Aut_{\F_q}(E)},
\]
where \(\Aut_{\F_q}(E)\) is the full automorphism group of the
unpointed genus one curve.  For genus two, collect the counts at each
number \(N\) of rational points of the plane quartic as
\[
 \mathsf M_{2,q}(N)
 =\sum_{\tau\in\{\mathrm{split},\mathrm{cusp},\mathrm{nonsplit}\}}
   \ \sum_{\substack{[C]\in\mathcal M_2(\F_q)\\
               \#C(\F_q)-\Delta_\tau=N}}
   \frac{m_\tau(C)}{\#\Aut_{\F_q}(C)}.
\]
Here the multiplicities and shifts are those of
Proposition~\ref{prop:g2-fixed-counts} and the preceding table.

\begin{thm}\label{thm:full-weight-enumerator}
The Hamming weight enumerator is
\begin{align*}
 W_{\cC_{2,4}}(X,Y)
={}&X^{q^2+q+1}
 +(q-1)\sum_z R_z(q)X^zY^{q^2+q+1-z}
\\*
&+(q^2+q+1)X^{q+1}
 \bigl(
  W_3^{A,\mathrm{cusp}}(X,Y)+W_3^{A,\mathrm{split}}(X,Y)\\*
&\hspace{3.6cm}{}+W_3^{A,\mathrm{nonsplit}}(X,Y)
  +W_3^{A,\mathrm{sm}}(X,Y)
 \bigr)
\\
&+(q-1)\#\PGL_3(\F_q)\sum_{N=0}^{q^2+q+1}
 \mathsf M_{2,q}(N)X^NY^{q^2+q+1-N}
\\
&+(q-1)\#\PGL_3(\F_q)\sum_{t^2\le4q}\sum_{d=-2}^{2}
 \mathsf B_{q,d}(t)
 X^{q+1-t-d}Y^{q^2+t+d}
\\
&+\frac{(q-1)\#\PGL_3(\F_q)}{\#\PGL_2(\F_q)}
 \sum_{d=-3}^{3}C_d(q)X^{q+1-d}Y^{q^2+d}
\\
&+(q-1)\#\PGL_3(\F_q)\sum_{n=0}^{q^2+q+1}
 S_q(n)X^nY^{q^2+q+1-n},
\end{align*}
where \(z\) runs through Table~\ref{tab:factorable-collected},
\(C_d(q)\) is given by Proposition~\ref{prop:g0-trace-shift-counts},
and \(S_q(n)\) is the smooth genus three distribution defined above.
\end{thm}

\begin{proof}
The zero form gives the first term, and
\eqref{eq:factorable-enumerator} gives the two reducible terms.
Grouping \eqref{eq:g2-weight-contribution} by \(N\) and
\eqref{eq:g1-weight-enumerator} by the normalization trace gives
the genus two and genus one terms.  Equation~\eqref{eq:g0-seven-term-weight}
gives the genus zero term, and the canonical plane models counted by
\(S_q(n)\) give the smooth term.
\end{proof}

At \(q=5,7\), the accompanying code evaluates every contribution,
giving complete weight enumerators.

\Needspace{6\baselineskip}
In particular, setting \(C_j(q)=0\) for \(|j|>3\), the singular
absolutely irreducible contribution has coefficients
\begin{align*}
 [X^NY^{q^2+q+1-N}]W_{\mathrm{irr,sing}}(X,Y)
 ={}&(q-1)\#\PGL_3(\F_q)\Bigg(
 \mathsf M_{2,q}(N)
 \\*
 &\quad+\sum_{d=-2}^{2}
 \mathsf B_{q,d}(q+1-N-d)
 +\frac{C_{q+1-N}(q)}{\#\PGL_2(\F_q)}\Bigg).
\end{align*}

The incidence coefficients and genus zero counts are explicit.  The next
two subsections evaluate \(\mathsf B_{q,d}(t)\) and
\(\mathsf M_{2,q}(N)\).  For general \(q\), the genus three
distribution remains in the form \(S_q(n)\).

\subsection{Genus one: formulas involving class numbers}

For a negative quadratic discriminant \(D\), let \(h(D)\) be the class
number of the order of discriminant \(D\), and set
\[
 h_w(D)=
 \begin{cases}
  h(D)/3,&D=-3,\\
  h(D)/2,&D=-4,\\
 h(D),&D<0,\ D\equiv0,1\pmod4,\ D\ne-3,-4,\\
 0,&\text{otherwise},
 \end{cases}
\]
The Hurwitz--Kronecker class number in this convention is
\(\mathscr H(D)=\sum_{\substack{a\ge1\\a^2\mid D}}h_w(D/a^2)\).
For elliptic curves, \(\Aut_{\F_q}(E,O)\) denotes the usual
automorphism group fixing the origin.  This differs from the full curve
automorphism group used in \(\mathsf B_{q,d}(t)\), which includes
rational translations.  Define the weighted count of elliptic curves
with trace \(t\) by
\[
 \mathscr E_q(t)
 =\sum_{\substack{[E,O]\;\mathrm{elliptic}/\F_q\\t_E=t}}
   \frac1{\#\Aut_{\F_q}(E,O)}.
\]
The corresponding weighted count for curves with full rational
two torsion is
\[
 \mathscr F_q(t)
 =\sum_{\substack{[E,O]\;\mathrm{elliptic}/\F_q\\t_E=t,\;
       E[2](\F_q)\simeq(\mathbf Z/2\mathbf Z)^2}}
   \frac1{\#\Aut_{\F_q}(E,O)}.
\]
Schoof gives the corresponding unweighted counts by trace and rational
torsion \cite[Theorem~4.6 and Lemma~4.8]{schoof1987nonsingular}.
We use the weighted formulas of Kaplan--Petrow recalled below.

\begin{lemm}[Kaplan--Petrow {\cite[\S3, Lemmas~2--3]{kaplanpetrow2017trace}}]
\label{lem:elliptic-class-number-counts}
For every integer \(t\), the weighted elliptic curve counts are
\begin{equation*}
\mathscr E_q(t)=
\begin{cases}
 \frac12\mathscr H(t^2-4q),
   &t^2<4q,\ p\nmid t,\\
 \frac12\mathscr H(-4p),
   &r\ \text{odd},\ t=0,\\
 \frac14\left(1-\left(\frac{-4}{p}\right)\right),
   &r\ \text{even},\ t=0,\\
 \frac16\left(1-\left(\frac{-3}{p}\right)\right),
   &r\ \text{even},\ t^2=q,\\
 \frac{p-1}{24},
   &r\ \text{even},\ t^2=4q,\\
 0,&\text{otherwise},
\end{cases}
\end{equation*}
and
\begin{equation*}
\mathscr F_q(t)=
\begin{cases}
 \frac12\mathscr H((t^2-4q)/4),
  &t^2<4q,\ p\nmid t,\ t\equiv q+1\pmod4,\\
 \frac12h_w(-p),
  &q\equiv3\pmod4,\ t=0,\\
 \mathscr E_q(t),&t^2=4q,\\
 0,&\text{otherwise}.
\end{cases}
\end{equation*}
Here \(\left(\frac{\cdot}{p}\right)\) is the Legendre symbol.
\end{lemm}

\begin{proof}
Kaplan--Petrow assign to \([E,O]\) the probability
\((q\#\Aut_{\F_q}(E,O))^{-1}\).  Multiplying their formulas by
\(q\) gives the two identities in our normalization.
\end{proof}

\begin{prop}\label{prop:g1-weighted-trace-counts}
For \(-2\le d\le2\), the coefficients in
Theorem~\ref{thm:full-weight-enumerator} are

\begin{align*}
\mathsf B_{q,2}(t)
 &=\frac{q^3-3q^2t-4q^2+3qt^2+8qt+4q-t^3-4t^2-4t
 -\mathbf 1_{2\mid t}}8\mathscr E_q(t)-\frac14\mathscr F_q(t)\\
\mathsf B_{q,1}(t)
 &=\frac{2q^2-3qt-4q+t^2+4t+\mathbf 1_{2\mid t}}2
 \mathscr E_q(t)+\mathscr F_q(t)\\
\mathsf B_{q,0}(t)
 &=\frac{3q^3+q^2t-2q^2-qt^2+6q+t^3+2t^2-2t-4
 -3\cdot\mathbf 1_{2\mid t}}4\mathscr E_q(t)-\frac32\mathscr F_q(t)\\
\mathsf B_{q,-1}(t)
 &=\frac{2q^2+qt-t^2+\mathbf 1_{2\mid t}}2
 \mathscr E_q(t)+\mathscr F_q(t)\\
\mathsf B_{q,-2}(t)
 &=\frac{q^3+q^2t-qt^2-t^3-\mathbf 1_{2\mid t}}8
 \mathscr E_q(t)-\frac14\mathscr F_q(t).
\end{align*}
\end{prop}

\Needspace{3\baselineskip}
\begin{proof}
We obtain \(\mathsf B_{q,d}(t)\) by weighting the fixed-curve counts of
Proposition~\ref{prop:g1-fixed-counts} by the full automorphism group.
For a curve of trace \(t\), each automorphism of the underlying curve
is uniquely a rational translation followed by an automorphism fixing
\(O\), so
\(\#\Aut_{\F_q}(E)=(q+1-t)\#\Aut_{\F_q}(E,O)\).
Translations identify all choices of rational origin, so forgetting
\(O\) gives a bijection between pointed elliptic curve classes and
unpointed genus one curve classes.  The pointed automorphism weight is \(q+1-t\) times the unpointed
weight, however.  Consequently,
\[
 \mathsf B_{q,d}(t)
 =\sum_{\substack{[E,O]\\t_E=t}}
   \frac{B_d(E)/(q+1-t)}{\#\Aut_{\F_q}(E,O)}.
\]
The factor \(\#E(\F_q)=q+1-t\) in every formula of
Proposition~\ref{prop:g1-fixed-counts} has now canceled.  The quotient is
linear in \(\#E[2](\F_q)\), with coefficients depending only on \(q,t\).

We therefore need the weighted sums of \(1\) and \(\#E[2](\F_q)\).
The first is \(\mathscr E_q(t)\).  For the second, observe that
\(\#E(\F_q)=q+1-t\) is odd when \(t\) is odd, so only the identity is
rational two torsion.  When \(t\) is even, the group has even order,
and \(\#E[2](\F_q)\) is either two or four.  It is four precisely for
the curves counted by \(\mathscr F_q(t)\).  At even trace, these curves
therefore contribute two more torsion points than the other curves.
This gives
\[
 \sum_{\substack{[E,O]\\t_E=t}}
 \frac{\#E[2](\F_q)}{\#\Aut_{\F_q}(E,O)}
 =(1+\mathbf 1_{2\mid t})\mathscr E_q(t)+2\mathscr F_q(t).
\]
In each quotient \(B_d(E)/(q+1-t)\), the part independent of
\(\#E[2](\F_q)\) thus multiplies \(\mathscr E_q(t)\), and the coefficient
of \(\#E[2](\F_q)\) multiplies the right-hand side above.

For example, when \(d=2\), Proposition~\ref{prop:g1-fixed-counts} gives
\(B_2(E)/(q+1-t)\) as
\(\bigl((q+1-t)^3-7(q+1-t)^2+15(q+1-t)-8\bigr)/8\)
minus \(\#E[2](\F_q)/8\).
The polynomial term multiplies \(\mathscr E_q(t)\), while the torsion
term contributes
\(-\frac{1+\mathbf 1_{2\mid t}}8\mathscr E_q(t)
 -\frac14\mathscr F_q(t)\).
Adding and expanding gives the stated formula for \(\mathsf B_{q,2}(t)\).
Applying the same rule to \(B_1(E),B_{-1}(E),B_{-2}(E)\) gives the other
three nonzero shifts.  For shift zero, the total count in
Proposition~\ref{prop:g1-fixed-counts}, after division by \(q+1-t\), is
\(q^3+q^2+t-1\).  It is independent of the two torsion, so
\(\sum_{d=-2}^{2}\mathsf B_{q,d}(t)
=(q^3+q^2+t-1)\mathscr E_q(t)\).
Subtracting the four nonzero shifts gives \(\mathsf B_{q,0}(t)\).
These substitutions and the class number evaluation are implemented in
the accompanying code.
\end{proof}

The same class number counts evaluate the smooth affine cubic term:
grouping \eqref{eq:affine-cubic-smooth} by trace gives
\[
 W_3^{A,\mathrm{sm}}(X,Y)
 =q^3(q-1)^2(q^2-1)
  \sum_{t^2\le4q}\mathscr E_q(t)W^\alpha(X,Y,t).
\]

\subsection{Genus two: sums over hyperelliptic models}

Define the three weighted distributions
\begin{align*}
 \mathsf H_q(t)
 &=\sum_{\substack{[C]\in\mathcal M_2(\F_q)\\q+1-\#C(\F_q)=t}}
   \frac1{\#\Aut_{\F_q}(C)},\\
 \mathsf V_q(t)
 &=\sum_{\substack{[C]\in\mathcal M_2(\F_q)\\q+1-\#C(\F_q)=t}}
   \frac{\#C(\F_{q^2})}{\#\Aut_{\F_q}(C)},\\
 \mathsf W_q(t)
 &=\sum_{\substack{[C]\in\mathcal M_2(\F_q)\\q+1-\#C(\F_q)=t}}
   \frac{\#\operatorname{Weier}(C)(\F_q)}{\#\Aut_{\F_q}(C)}.
\end{align*}
Each of these distributions vanishes when \(|t|>4\sqrt q\).

\begin{prop}\label{prop:g2-weighted-point-counts}
The genus two coefficient in Theorem~\ref{thm:full-weight-enumerator} is
\begin{align*}
 \mathsf M_{2,q}(N)
={}&\frac{(N^2-1)\mathsf H_q(q-N)+\mathsf W_q(q-N)}2
 \\
 &+N\mathsf H_q(q+1-N)-\mathsf W_q(q+1-N)
 \\
 &+\frac{\mathsf V_q(q+2-N)+\mathsf W_q(q+2-N)
            -(2q+2)\mathsf H_q(q+2-N)}2.
\end{align*}
\end{prop}

\begin{proof}
A split node identifies two rational points of the normalization with
one rational point of the quartic.  Thus a quartic with \(N\) rational
points has \(\#C(\F_q)=N+1\), so the normalization trace is \(q-N\).
Proposition~\ref{prop:g2-fixed-counts} gives
\(m_{\mathrm{split}}(C)
=(N^2-1+\#\operatorname{Weier}(C)(\F_q))/2\).
Summing with weight \(1/\#\Aut_{\F_q}(C)\) replaces the constant
term by \((N^2-1)\mathsf H_q(q-N)/2\) and the Weierstrass term by
\(\mathsf W_q(q-N)/2\), giving the first term of the formula.

A cusp has one rational preimage, so \(\#C(\F_q)=N\) and the trace
is \(q+1-N\).  The count
\(m_{\mathrm{cusp}}(C)=N-\#\operatorname{Weier}(C)(\F_q)\)
therefore gives
\(N\mathsf H_q(q+1-N)-\mathsf W_q(q+1-N)\).

A nonsplit node is rational but has no rational preimage, so
\(\#C(\F_q)=N-1\) and the trace is \(q+2-N\).
The nonsplit count is
\(m_{\mathrm{nonsplit}}(C)
=(\#C(\F_{q^2})+\#\operatorname{Weier}(C)(\F_q)-2q-2)/2\).
The three summands yield respectively \(\mathsf V_q/2\),
\(\mathsf W_q/2\), and \(-(q+1)\mathsf H_q\), all evaluated at
\(q+2-N\), this is the third term.
\end{proof}

To evaluate these distributions, let \(\mathcal P_2(q)\) be the set of
all squarefree polynomials in \(\F_q[x]\) of degree five or six.
The polynomials are not restricted to be monic.
For \(f\in\mathcal P_2(q)\), let \(C_f\) be the smooth projective
model of \(y^2=f(x)\).

\begin{lemm}[{\cite[\S3, equation~(3.2)]{bergstrom2009hyperelliptic}}]
\label{lem:genus-two-model-count}
For every rational-valued isomorphism invariant function \(A\) on
smooth genus two curves over \(\F_q\),
\[
 \sum_{[C]\in\mathcal M_2(\F_q)}
 \frac{A(C)}{\#\Aut_{\F_q}(C)}
 =\frac1{\#\GL_2(\F_q)}\sum_{f\in\mathcal P_2(q)}A(C_f).
\]
\end{lemm}

Applying Lemma~\ref{lem:genus-two-model-count} to the weights defining
\(\mathsf H_q(t),\mathsf V_q(t),\mathsf W_q(t)\) expresses the
coefficient in Proposition~\ref{prop:g2-weighted-point-counts} as a
finite sum over \(\mathcal P_2(q)\).
The regrouping and the model sums at \(q=5,7\) are implemented in the
accompanying code.

For \(q=5,7\), the accompanying code also evaluates \(S_q(n)\).
Curves with a rational point are counted using normal forms obtained by
fixing the point and its tangent.  The normal forms and counting weights
are justified in the accompanying code.
The terms \(S_q(0)\) use the classification of Faber, Grantham,
and Howe and its accompanying data
\cite[\S8]{faber-grantham-howe2025gonality}.
The resulting complete enumerators satisfy
\(W_{\cC_{2,4}}(1,1)=q^{15}\).

\section*{Acknowledgments}
The author thanks Nathan Kaplan for his guidance, helpful discussions, and
comments on the manuscript.  The author received support from NSF grant
DMS 2154223.

\enlargethispage{\baselineskip}
\paragraph{AI acknowledgment.}
During the preparation of this manuscript, the author used OpenAI's
ChatGPT and Codex, with the \mbox{GPT-5.6 Sol}, \mbox{GPT-6 Astra}, and
\mbox{GPT-6.1 Sol} models, primarily to discuss the results and arguments
for singular quartic curves and to write SageMath code for the counting
calculations.
The author is responsible for the mathematical arguments, computations,
and conclusions presented in this paper.

\appendix
\section{Incidence tables for reducible quartic forms}
\label{app:reducible-tables}
\setcounter{table}{0}
\renewcommand{\thetable}{\thesection.\arabic{table}}
\renewcommand{\theHtable}{appendix.\arabic{table}}

The tables below list the counts for each case in
\Cref{sec:reducible}.  The geometric arguments are given in the main text.

\FourLineCountingTable

\DoubleLineConicCountingTable

\TwoLinesConicCountingTable

\TwoConicsCountingTable

\ConjugateLinesConicCountingTable

\ConjugateConicsCountingTable

\begingroup
\setlength{\bibitemsep}{3pt plus 1pt minus 1pt}
\printbibliography
\endgroup

\end{document}